\documentclass[11pt, reqno]{amsart}
\usepackage[hmargin=35mm,vmargin=35mm]{geometry}
\usepackage{amsmath, amsthm, amssymb, mathabx, verbatim, setspace, enumerate,mathtools}
\usepackage[mathscr]{euscript}
\usepackage[all,cmtip]{xy}
\usepackage{etoolbox}
\usepackage{xcolor}
\usepackage{tikz-cd}
\usepackage{hyperref}

\numberwithin{equation}{section}

\newcommand{\ds}{\displaystyle}

\newcommand{\CC}{\mathbb{C}}

\newcommand{\ZZ}{\mathbb{Z}}
\newcommand{\RR}{\mathbb{R}}

\newcommand{\F}{\mathcal{F}}

\newcommand{\HH}{\mathcal{H}}

\newcommand{\C}{\mathcal{C}}

\newcommand{\eps}{\varepsilon}

\newtheorem{Thm}{Theorem}[section]
\newtheorem{Prop}[Thm]{Proposition}
\newtheorem{Lem}[Thm]{Lemma}
\newtheorem{Cor}[Thm]{Corollary}

\theoremstyle{definition}
\newtheorem{rmk}[Thm]{Remark}

\DeclareMathOperator{\SL}{SL}

\usepackage[OT2,T1]{fontenc}
\DeclareSymbolFont{cyrletters}{OT2}{wncyr}{m}{n}
\DeclareMathSymbol{\Sha}{\mathalpha}{cyrletters}{"58}

\title{Binary forms with covariant points close to the real axis}

\author{Eugenia Rosu}
\date{}

\begin{document}
\maketitle
\begin{abstract}  For a real binary form $F(X, Z)$, Stoll and Cremona have defined a reduction theory using the action of the modular group $\SL_2(\ZZ)$, and associated to each binary form a covariant point $z(F)$ located in the upper half plane. When the point $z(F)$ is close to the real axis, then at least half of the roots will be on a circle of small radius $r$. Conversely, we find conditions depending on the radius $r$ such that the covariant point $z(F)$ to be close to the real axis.  The results have further applications to improving the reduction algorithm for binary forms of Stoll and Cremona.
\end{abstract}

\section{Introduction}
Let $F(X, Z)=a_{0}\left(X-\alpha_{1} Z\right) \ldots\left(X-\alpha_{n} Z\right)$ be a real binary form of degree $n$ with roots $\alpha_{1}, \ldots, \alpha_{n}\in \CC$. In \cite{SC}, Stoll and Cremona have defined a reduction theory for complex and real binary forms, inspired by the work of Julia \cite{J} for $n=3, 4$. Their reduction theory has as motivation finding equivalent binary forms with smaller coefficients and, over $\ZZ$, finding forms with the same discriminant.  This has further applications to hyperelliptic curves.

 Stoll and Cremona defined an equivalence of binary forms by letting the modular group $\SL_2(\ZZ)$ act by linear transformations on the variables:
\[
(F g)(X, Z)=F(aX+bZ, cX+dZ), \text{ for } g\in \SL_2(\ZZ).
\]
 To each real form $F$ one can associate a covariant point $z(F)=(t, u)$ in the upper half plane $\HH$, where by covariance we understand that 
\[
z(F g^{-1})=g(z(F)), 
\]
for all $g \in \SL_2(\ZZ)$. Here $\SL_2(\ZZ)$ acts on the left on $\HH$ by fractional linear transformations. For each $\SL_2(\ZZ)$-orbit in $\HH$ there is a representative in the standard fundamental domain 
\[
\mathfrak{F}=\left\{z \in \HH:|z| \geq 1,-\frac{1}{2} \leq \operatorname{Re}(z) \leq \frac{1}{2}\right\}.
\] The representative in $\mathfrak{F}$ is unique, except when it is on the boundary of $\mathfrak{F}$ in which case we get at most two representatives. We say that $F$ is {\it reduced} iff $z(F) \in \mathfrak{F}$.

Stoll and Cremona give in \cite{SC} a beautiful geometric interpretation of the covariant $z(F)=(t, u)$: it is the unique point in the hyperbolic half-space $\HH_3$ such that the sum of its hyperbolic distances from all the roots of $F$ is minimal. Moreover, they proved that $z(F)$ is uniquely determined by the roots of the form $F$ through a system of equations (see \eqref{eqn:1.1} and \eqref{eqn:1.2}) and provided an algorithm for finding a covariant point associated to a reduced binary form equivalent to $F$.  

\smallskip
In the current paper we are interested in the case of the covariant point $z(F)=(t, u)\in \HH$ being close to the real axis (i.e.  $u$ is small). We find equivalent conditions depending on disks containing the roots of the binary form $F(X, Z)$. 

First, when $z(F)=(t, u)$ is near the real axis, we show that at least half of the roots are on a disk of radius considerably small. Moreover, the reverse is also true if more than half of the roots on a disk of small radius:

\begin{Prop}\label{thm:intro_1}
\begin{enumerate}[(i)]
	\item If $u\leq\epsilon$, then at least half of the roots are on a disk of radius $\ds\epsilon \sqrt{n}$ centered at $t$.

	\item If more than half of the roots are on a disk of radius $r\leq \epsilon$, then $u\leq 2\epsilon\sqrt{n}$.
\end{enumerate}
\end{Prop}

\medskip

However, when precisely half of the roots are on a circle of small radius $r_1$ (with a center $c_1$), it is no longer enough for $r_1$ to be small in order for $u$ to also be small. Instead, we have to construct the disk of minimum radius $r_2$, with a center $c_2$, containing the other half of the roots, and add one extra condition:

\begin{itemize}
	\item when the distance between the two centers is large enough, $u$ being small is equivalent to the ratio $r_1/r_2$ being small as well.
	\item when the two centers are closer together, $u$ being small is related to $r_1r_2$ being small.

\end{itemize}

First, if the centers are not too close to each other:

\begin{Thm}\label{thm:intro_2} Let $\ds \epsilon<\frac{1}{100n^2(2n+3)}$. Assume that $r_1\leq \epsilon$ and
\[
|c_1-c_2|\geq \max(4\sqrt{\epsilon}, 4r_2).
\]
Then we have:
  
  (i) $u\leq\epsilon$ implies $\ds\frac{r_1}{r_2}\leq\frac{10n\epsilon}{3}$.
  
  (ii) $\ds \frac{r_1}{r_2}\leq\frac{10n\epsilon}{3}$ implies $\ds u\leq \frac{20n\epsilon}{7}$.

\end{Thm}

Better bounds are presented in Propositions \ref{prop:4.9} and \ref{prop:4.10}, as well as Proposition \ref{prop:r2small_u} and Proposition \ref{prop:r2notsmall_u}. 

\smallskip
We note that the case of $r_2\leq \sqrt{\epsilon}$ and $|c_1-c_2|\leq 4\sqrt{\epsilon}$ reduces to the case of all the roots on a circle of small radius $3\sqrt{\epsilon}$, which implies from Proposition \ref{thm:intro_1}(ii) that $u\leq 6\sqrt{n\epsilon}$. 

We are left only with the case of the centers being closer together and $r_2\geq \sqrt{\epsilon}$ (thus not too small). In this case the condition of $u$ being small is related to $r_1r_2$ being bounded in almost all cases. The only case left out is that of the center $c_1$ almost being on the second circle (i.e. $|c_1-c_2|\sim r_2$) or $\ds |t-c_1|\leq r_1+\frac{u}{\sqrt{n}}$.

\begin{Thm}\label{thm:intro_3} Let $\ds\epsilon<\frac{1}{100n^2(2n+3)}$. Assume that:
\begin{itemize}
	\item $r_1\leq \epsilon$, $r_2\geq \sqrt{\epsilon}$,
	\item $|c_1-c_2|\leq 2r_2+r_2\sqrt{\epsilon}+r_1$. \end{itemize}

(1) Assume further
\[
\ds |t-c_1|\geq r_1+\frac{u}{\sqrt{n}}, \text{ and } \ds \left |\frac{|c_1-c_2|}{r_2}-1\right|\geq 64n\sqrt{\epsilon}.
\] 

If $u\leq \epsilon$, then we have the bound $\ds r_1r_2<\frac{1}{64n^2}$.

\smallskip
(2) If $r_1r_2<\epsilon^2$, then $u<2\sqrt{\epsilon}$.

\end{Thm} 

The interest of the paper is mainly theoretic, but the question is motivated by the following computational question. In order to find a reduced form that is $\SL_2(\ZZ)$-equivalent to $F(X, Z)$, Stoll and Cremona presented the following algorithm. First calculate $z:=z(F)$ and until $z \in \mathfrak{F}$, do:
\begin{enumerate}
\item take $m=\text{round}(t)$, the integer closest to $t$; set $F(X, Z):=F(X+m Z, Z)$, which implies $z:=z-m$;

\item if $|z|<1$, set $F(X, Z):=F(Z,-X)$, which implies $z:=-1/z$.
\end{enumerate}
\medskip

However, in the case when $z(F)$ is very close to the real axis, it is hard to compute $z(F)$ with reasonable accuracy. In the present paper we are improving the algorithm in this case by finding an approximation for the covariant point. The goal is to take $m=\text{round}(\text{Re}(c))$ instead of $m=\text{round}(t)$ in the algorithm, where $c$ is the center of the circle with at least half of the roots. The advantage is that the coordinates of $c$ can be easily computed (optimally using the algorithm of Megiddo). We show:


\begin{Thm}\label{thm:intro_4}
Let $F(X, Z)$ be a real binary form with precisely $k\geq n/2$ roots on a disk centered at $c$ and small radius $r\leq \epsilon$. If $u\leq \epsilon$, then:

\begin{enumerate}
	\item If $k>n/2$, then $|t-c| \leq(2n+3) \epsilon$.
	\item If $k=n/2$, then $\ds 
|t-c| \leq \frac{1}{2\sqrt{n}}+ \epsilon.$
\end{enumerate}
\end{Thm} 

This implies that $t$ and $\text{Re}(c)$ have the same closest integer, thus we can use $\text{Re}(c)$ instead of the $t$ coordinate in the algorithm of Stoll and Cremona \cite{SC}. This is described in detail in Section \ref{sec:algorithm}.

\medskip

The paper is structured as follows. In Section \ref{sec:2} we show that when $u$ is small, we must have at least half of the roots on a circle of small radius. In Section \ref{sec:3}, assuming there are more than $n/2$ roots in a circle of small radius $r$ and center $c$, we find a bound for $u$ with respect to $r$ and show that $t$ and $c$ are close. 

In the remaining sections we consider the case of precisely $n/2$ roots on a circle of small radius $r_1$, centered at $c_1$.
In Section \ref{sec:4} we find bounds for $u$ and the ratio of the two radii $r_1/r_2$. In Section \ref{sec:5} we assume that $u$ is small, and show that depending on the distance of the two centers, either the ratio $r_1/r_2$ is also small, or the product $r_1r_2$ is bounded. Moreover, we show that $t$ and $c_1$ will have the same closest integer. In Section \ref{sec:6} we show the converse statements. Under the assumption that either the ratio $r_1/r_2$ or the product $r_1r_2$ is small, we show that $u$ will also be small and $t$ will be close to $c_1$. In Section \ref{sec:algorithm} we discuss the applications to the algorithm of Stoll and Cremona from \cite{SC}.

\subsection*{Acknowledgements}
A portion of the paper has been completed in $2009$ during a two month-internship at the Bayreuth University under the supervision of Michael Stoll. I would like to thank Michael Stoll for suggesting the project and for his comments.

\section{The number of roots in a disk centered at $t$}\label{sec:2}

Let $F(X, Z)=a_0(X-\alpha_1Z)\dots (X-\alpha_nZ)$ be a real binary form of degree $n$. In \cite{SC}, Stoll and Cremona have associated to $F$ a covariant point $z(F)=(t, u)$ in the upper half plane $\HH=\{x+yi\in \CC: y>0\}$, which is uniquely determined by the system of equations:

\begin{equation}\label{eqn:1.1}
 \sum_{i=1}^{n} \frac{u^{2}}{|t-\alpha_{i}|^{2}+u^{2}}=\frac{n}{2},
\end{equation}
\begin{equation}\label{eqn:1.2}
 \sum_{i=1}^{n} \frac{t-\alpha_{i}}{|t-\alpha_{i}|^{2}+u^{2}}=0.
\end{equation}

Moreover, they present a geometric description of the point $z(F)$. Note that the upper half plane $\HH$ can be considered as a vertical cross-section of the hyperbolic upper half-space 
\[
\HH_3=\left\{(z, u): z \in \mathbb{C}, u \in \RR_{>0}\right\}
\]
 by identifying $(t, u) \in \HH_3$ with $t+i u \in \HH$. We can view the roots of $F$ as lying on the boundary of $\HH_3$ identified with $\mathbb{P}^{1}(\mathbb{C})$; take the hyperbolic geodesics through each root and $z(F)$. Then the point $z(F)$ is characterised by the property that the unit tangent vectors at $z(F)$ in the directions of the roots of $F$ add up to zero.

\medskip

In this section, we show that when $u$ is close to $0$, there are at least $n/2$ roots on a disk of small radius. We do this by looking at the number of roots in a circle centered at $t$.
 
\begin{Prop}\label{prop:2.1} Denote by $k$ the number of roots of the form $F(X, Z)$ inside the disk of center $t$ and radius $R$. Then we can bound $k$ in the following way:

\begin{equation*}\label{eqn:2.1}
\frac{n}{2} \cdot\left(1-\frac{u^{2}}{R^{2}}\right)
<k
 \leq \frac{n}{2}\left(1+\frac{R^{2}}{u^{2}}\right)
\end{equation*}
\end{Prop}

\noindent{\bf Proof:} Let $\alpha_{1}, \ldots, \alpha_{k}$ be the $k$ roots that are on the disk of center $t$ and radius $R$.

Then $|t-\alpha_{j}| \leq R$ for $1 \leq j \leq k$ which implies $\ds \frac{u^{2}}{u^{2}+|t-\alpha_{j}|^{2}} \geq \frac{1}{1+R^{2} / u^{2}}$. Using the obvious bound $\ds\frac{u^{2}}{u^{2}+|t-\alpha_{j}|^{2}}>0$ for $k+1 \leq j \leq n$ and the equation \eqref{eqn:1.1} we get

\[
\frac{n}{2}=\sum_{j=1}^{k} \frac{u^{2}}{u^{2}+|t-\alpha_{j}|^{2}}+\sum_{j=k+1}^{n} \frac{u^{2}}{u^{2}+|t-\alpha_{j}|^{2}} \geq \frac{k}{1+R^{2} / u^{2}},
\]
which gives the right-hand side of the inequality. Moreover, note that the inequality is strict if $k<n$, and the equality occurs only for $R=u$.
\medskip

For $k+1 \leq j \leq n$ we have $|t-\alpha_{j}|>R$, implying $\ds \frac{u^{2}}{u^{2}+|t-\alpha_{j}|^{2}}<\frac{1}{1+R^{2} / u^{2}}$, while clearly $\ds \frac{u^{2}}{u^{2}+|t-\alpha_{j}|^{2}} \leq 1$ for $1 \leq j \leq k$. Then using \eqref{eqn:1.1} we obtain:
\[
\frac{n}{2}=\sum_{j=1}^{k} \frac{u^{2}}{u^{2}+|t-\alpha_{j}|^{2}}+\sum_{j=k+1}^{n} \frac{u^{2}}{u^{2}+|t-\alpha_{j}|^{2}} < k+\frac{n-k}{1+R^{2} / u^{2}}
\]
Easy calculations imply the left-hand side of the inequality, and it is easy to see that the equality is not possible. Note that the same inequalities are true for complex forms.

\medskip

By plugging in $R=u\sqrt{n}$, respectively $R=u/\sqrt{n}$, we get immediately:

\begin{Cor}\label{cor:2.2}

\begin{enumerate} 
\item There are at least $n/2$ roots on the disk of center $t$ and radius $u \sqrt{n}$.

\item There are at most $n/2$ roots on the disk of center $t$ and radius $u/\sqrt{n}$.
\end{enumerate}
\end{Cor}
\noindent{\bf Proof:} Part (1) follows by taking $R=u \sqrt{n}$ in the left-hand side of \eqref{eqn:2.1}.

For (2), take $R=u / \sqrt{n}$, and as $R \neq u$, the right-hand side of the inequality \eqref{eqn:2.1} is strict and we get $k<\frac{n+1}{2}$.

\bigskip

For $u<\epsilon$, Corollary \ref{cor:2.2} (1) implies immediately:

\begin{Cor}\label{cor:usmall} If $u<\epsilon$, then at least half of the roots are on a disk of radius $\ds\epsilon \sqrt{n}$ centered at $t$.
\end{Cor}

This means that when $u \approx 0$, then at least half of the roots are on a disk of small radius centered at a real number.

Also, in Section \ref{sec:4} we will be interested in the disk of minimum radius that contains exactly half of the roots. The following result which is an immediate consequence of Corollary \ref{cor:2.2} (2) will be useful.

\begin{Cor} Let $r$ be the radius of the smallest disk containing $n/2$ of the roots. Then $r \leq u \sqrt{n}$.
\end{Cor}

\smallskip

\section{More than $n/2$ roots on a circle of small radius}\label{sec:3}

Let $F(X, Z)=a_{0}\left(X-\alpha_{1} Z\right) \ldots\left(X-\alpha_{n} Z\right)$ be a real binary form such that exactly $k>n/2$ of the roots, say $\alpha_{1}, \alpha_{2}, \ldots, \alpha_{k}$, are on a disk of center $c$, radius $r$. We let $z(F)=(t, u)\in \HH_3$ be the associated covariant point.

 In this section we will prove that when more than half of the roots are on a disk of radius $r$ and center $c$, then $u$ is small and the coordinate $t$ and center $c$ are close to each other.
 
\subsection{Bound for $u$}\label{sec:3.1}

 The goal of this section is to prove the following bound for $u$:

\begin{Prop}\label{prop:3.1} If there are $k>n/2$ roots on the disk of center $c$ and radius $r$, then we have the upper bound $u \leq 2 r \sqrt{n}.$
\end{Prop}

We note Proposition \ref{prop:3.1} gives Proposition \ref{thm:intro_1}(ii), while Corollary \ref{cor:usmall} gives Proposition \ref{thm:intro_1}(i).

The statement of Proposition \ref{prop:3.1} is easy to see in the case of $|t-c| \leq r$, as there are more than $n/2$ roots on the disk of radius $2r$ and center $t$, and by Corollary \ref{cor:2.2} (2) we have
\begin{equation}\label{eqn:3.1}
\frac{u}{\sqrt{n}} \leq 2 r.
\end{equation}

In the the following we will look at the case $|t-c|>r$ and compute a bound for $u$. In order to do that, we define the binary form 
\[
G(X, Z)=a_{0}(X-\beta_{1} Z) \dots (X-\beta_{n}Z),
\] 
where $\ds \beta_{i}=\frac{t-\alpha_{i}}{u}$ for $1 \leq i \leq n$. The system of equations \eqref{eqn:1.1}, \eqref{eqn:1.2} becomes in terms of $\beta_{i}$ 's:
\begin{equation}\label{eqn:3.2}
 \sum_{i=1}^{n} \frac{1}{1+|\beta_{i}|^{2}}=\frac{n}{2},  
\end{equation}
\begin{equation}\label{eqn:3.3}
 \sum_{i=1}^{n} \frac{\beta_{i}}{1+|\beta_{i}|^{2}}=0.
\end{equation}

Note that the form $G(X, Z)=a_{0}\left(X-\beta_{1} Z\right) \ldots\left(X-\beta_{n} Z\right)$ has the unique covariant point $z(G)=(0,1)\in \HH_3$. We denote by $M_{j}=\left(m_{j}, n_{j}, p_{j}\right)$ the point on the unit sphere which is the inverse stereographic projection of $\beta_{j}, 1 \leq j \leq n$, and let $\vec{v}_j=\overrightarrow{OM}_{j}$ be the unit vector from the origin $O$ to $M_j$. Then explicitly $\ds \beta_{j}=\left(\frac{m_{j}}{1-p_{j}}, \frac{n_{j}}{1-p_{j}}, 0\right)$ for $1 \leq j \leq n$. 
\[
        \begin{tikzpicture}
            \coordinate (A) at (3,-0.25) node [yshift=-14mm, xshift=1mm, color=black] {Figure 1: stereographic projection};
            \coordinate (P) at (0,2);
            \coordinate (O) at (0,0);
;

            \draw (0:2cm)   arc[radius=2cm,start angle=0,end angle=180];
             
            \draw (180:2cm) arc[x radius=2cm, y radius=0.5cm, start angle=180,end angle=360];

            \draw [dashed] (180:2cm) 
            arc[start angle=180,delta angle=-180,x radius=2cm,y radius=0.5cm];

         	\draw [dashed] (150:2cm) coordinate(ul) -- (30:2cm) coordinate(ur);

            \draw (-4.5,-1) -- (3.5,-1) -- (4.5,1) node[anchor=south east] {\scriptsize} -- (ur) (ul) -- (-3.5,1) -- (-4.5,-1);

            \draw (A) -- (P) coordinate[pos=0.47](B);
     
            \path (A) node[circle, fill, inner sep=1pt, label=below:{\scriptsize$ \beta_j$}]{};
            \path (B) node[circle, fill, inner sep=1pt, label=left:{\scriptsize$M_j$}]{};
            \path (P) node[circle, fill, inner sep=1pt, label=above:{\scriptsize$ z(G)=(0,0,1) $}]{};
            
            \path (O) node[circle, fill, inner sep=1pt, label=above:{\scriptsize$ O$}]{};
            
            \draw[->, thick] (O) -- (B) node[inner sep=1pt, yshift=-9mm, xshift=-7mm, label=above:{\scriptsize$\vec{v}_j$}]{};
        \end{tikzpicture}
\]

We note that the tangent vector at $z(G)=(0, 0, 1)$ for the hyperbolic geodesic through $z(G)$ and $\beta_{j}$, in the direction of $\beta_j$, is the vertical translation of $\vec{v}_j$ by the unit vector from $O$ to $z(G)$. Stoll and Cremona have proved in \cite{SC} that in the upper half-space $\HH_3$ the sum of the unit tangent vectors at $z(G)$ in the directions of the roots of the binary form $G$ add up to zero, which is equivalent to
\[
\sum_{j=1}^{n} \vec{v}_{j}=\vec{0}.
\]

As the roots $\alpha_1, \dots, \alpha_k$ are on a disk of radius $r$ centered at $c$, the roots $\beta_{1}, \dots, \beta_{k}$ are on a disk of center $\ds c_0=\frac{t-c}{u}$ and radius $\ds r_{0}=\frac{r}{u}$. As we are in the case of $|c_0|>r_{0}$, all the tangent unit vectors at $z(G)=(0,1)$ are in the same half-space. For $1 \leq i, j \leq k$, denote by $\theta_{ij}$ the angle at $O$ between the vectors $\vec{v}_i$ and $\vec{v}_j$. The dot product of the two unit vectors $\vec{v}_i$ and $\vec{v}_j$ gives us $m_{i} m_{j}+n_{i} n_{j}+$ $p_{i} p_{j}=\cos \theta_{i j}$, and a simple calculation implies:
\begin{equation}\label{eqn:3.4}
|\beta_{i}-\beta_{j}|^{2}=\frac{2\left(1-\cos \theta_{i j}\right)}{\left(1-p_{i}\right)\left(1-p_{j}\right)}.
\end{equation}

\smallskip

We define:

\begin{itemize}
	\item $\theta=\max\limits_{1\leq i, j\leq k} \theta_{i j}$ the maximum angle between any two vectors $\vec{v}_i, \vec{v}_j, 1\leq i, j \leq k$,
	\item $\gamma=$ the greatest angle at $O$ on the smallest spherical cap containing all points $M_i$ for $1\leq i\leq k$.
\end{itemize}
 Note that $\theta, \gamma \leq \pi$, as all the points $M_{i}, 1 \leq i \leq k$, are in the same hemisphere. 
 
We will need the following two inequalities in order to find a bound for $u$. These inequalities will also be used in the proof of Proposition \ref{prop:4.7}.
\begin{Lem}\label{lem:3.2}
We have the inequalities:
\begin{enumerate}
	\item $\ds \sin{\frac{\theta}{2}} \geq \frac{\sqrt{3}}{2} \sin\frac{\gamma}{2}$
	\item $\ds |\sum\limits_{i=1}^{k} \vec{v}_i| \geq k \cos \frac{\gamma}{2}$.
\end{enumerate}
\end{Lem}
\noindent{\bf Proof:} (1) Denote by $C$ the center of the circle that determines the spherical cap. Then we can find two points on this circle, say $M_1, M_2$, such that the angle $\angle M_{1} C M_{2}$, is greater than or equal to $2 \pi/3$.

By looking at the two triangles $M_{1} O M_{2}$ and $M_{1} C M_{2}$, we can write the edge $M_1M_2$ as:
\[
M_{1} M_{2}=2 \sin \frac{\theta}{2}=2 CM_{1} \sin \frac{\angle M_{1} C M_{2}}{2}.
\]
 Also, by looking at the triangle 
$OCM_{1}$ we get $C M_{1}=\sin \frac{\gamma}{2}$. Together with $\sin \frac{\angle M_1CM_2}{2} \geq \sqrt{3}/2,$ then we get:
\[
2 \sin \frac{\gamma}{2} \sin \frac{\angle M_{1} C M_{2}}{2} \geq \sqrt{3} \sin \frac{\gamma}{2},
\]
which gives us the conclusion.

(2) To find the minimum norm $|\sum_{i=1}^k \vec{v_i}|$ it is enough to assume that all $M_1, \dots, M_k$ are on the circle that bounds the cap. Then we can write $\ds \sum_{i=1}^k\vec{v}_i=k \vec{w}+\sum_{i=1}^{k} \overrightarrow{M^{\prime} C_{i}}$, where $\vec{w}:=\overrightarrow{OC}$. As $\vec{w}$ has length $\ds \cos\frac{\gamma}{2}$, the minimum norm of $\ds \sum\limits_{i=1}^k\vec{v}_i$ is $\ds k \cos \frac{\gamma}{2}$, attained when the vertices $M_{1}, \dots, M_{k}$, form a regular $k$-gon centered at $C$.

\medskip

These inequalities help us find a stronger bound for $u$ in the case of $|t-c|>r$:
\begin{Prop}\label{prop:3.1.0} If there are $k>n/2$ roots on the disk of center $c$, radius $r$ and $|t-c|>r$, then we get the following upper bound for $u$:
\begin{equation}\label{eqn:3.6}
u \leq \frac{4 k r}{\sqrt{3 n(2 k-n)}}.
\end{equation}

\end{Prop}
\noindent{\bf Proof:} As $\sum\limits_{i=1}^{n} \vec{v}_i=\overrightarrow{0}$, we have the equality of norms $|\sum\limits_{i=1}^{k} \vec{v}_i|=|\sum\limits_{i=k+1}^{n} \vec{v}_i|$. As $\vec{v}_i$ are unit vectors, the maximum value of the norm $|\sum\limits_{i=k+1}^{n} \vec{v}_i|$ is $n-k$, and from Lemma \ref{lem:3.2}(2), the minimum value of the norm $|\sum\limits_{i=1}^{k} \vec{v}_i|$ is $k \cos \frac{\gamma}{2}$, implying that $k \cos \frac{\gamma}{2} \leq n-k$. This is equivalent to $\ds \sin ^{2} \frac{\gamma}{2} \geq \frac{(2k-n) n}{k^{2}}$. Combining this with Lemma \ref{lem:3.2}(1), we obtain
\begin{equation}\label{eqn:3.7}
1-\cos \theta=2 \sin ^{2} \frac{\theta}{2} \geq \frac{3}{2} \sin ^{2} \frac{\gamma}{2} \geq \frac{3 n(2 k-n)}{2 k^{2}}.
\end{equation}

On the other side, we have $|\beta_{i}-\beta_{j}| \leq 2r/u$ for any $i, j$, as $\beta_{i}$ and $\beta_{j}$ are on a disk of radius $r/u$.
The maximum angle $\theta$ between the tangent unit vectors on the unit sphere is attained for two vectors, say $\vec{v}_i$ and $\vec{v}_j$. Using \eqref{eqn:3.4} thus we obtain 
\[
\ds |\beta_{i}-\beta_{j}|^{2}=\frac{2(1-\cos \theta)}{\left(1-p_{i}\right)\left(1-p_{j}\right)} \leq \frac{4 r^{2}}{u^{2}}.
\]
Then using the inequality \eqref{eqn:3.7} and the fact that $1>p_{i}, p_{j} \geq-1$, we get:
\[
\frac{3 n(2 k-n)}{k^{2}} \leq 2(1-\cos \alpha) \leq\left(1-p_{1}\right)\left(1-p_{2}\right) \frac{4 r^{2}}{u^{2}} \leq 16 \frac{r^{2}}{u^{2}},
\]
and the conclusion follows.

\bigskip

Note that the inequality $\ds \frac{4kr}{\sqrt{3 n(2k-n)}} \leq 2 r \sqrt{n}$ together with \eqref{eqn:3.1} gives immediately Proposition \ref{prop:3.1}.

\subsection{Bound for $t$.}\label{sec:3.2}

As before, we will first find bounds using the binary real form $G(X, Z)$ with covariant point $z(G)=(0, 1)$, and use this to construct a bound for $t$.

\begin{Prop} If $G(X, Z)$ has precisely $k$ roots on a disk of center $c_0$ and radius $r_0$, then for $k>n/2$ we must have
\begin{equation}\label{eqn:3.8}
\frac{1}{\sqrt{n}}-r_{0} \leq|c_0| \leq \sqrt{n}+r_{0}.
\end{equation}

\end{Prop}
\noindent{\bf Proof:} Let $\beta_{1}, \dots, \beta_{k}$ be the roots on the disk centered at $c_0$ and radius $r_{0}$, which implies:
\begin{equation}\label{eqn:3.9}
|| c_0|-r_{0}| \leq|\beta_{i}| \leq |c_0|+r_{0}. 
\end{equation}

Using \eqref{eqn:3.2} and the right-hand side of \eqref{eqn:3.9} we obtain 
\[
\frac{n}{2}
=
\sum_{i=1}^{n} \frac{1}{1+|\beta_{i}|^{2}} \geq 
\sum_{i=1}^{k} \frac{1}{1+|\beta_{i}|^{2}}>
\frac{k}{1+\left(|c_0|+r_{0}\right)^{2}}.
\] 
After easy calculations and using the fact that $k \geq(n+1)/2$, we get $|c_0|+r_{0} \geq 1/\sqrt{n}$.

Similarly, using the left-hand side of \eqref{eqn:3.9} and the obvious bound $\ds\frac{1}{1+|\beta_{i}|^{2}}<1$ we obtain:
\[
\frac{n}{2}=\sum_{i=1}^{n} \frac{1}{1+|\beta_{i}|^{2}}<\frac{k}{1+\left(|c_0|-r_{0}\right)^{2}}+(n-k).
\]
The inequality is equivalent to $(k-n/2)\left(|c_0|-r_{0}\right)^{2} \leq n/2$, so $||c_0|-r_{0}|^{2} \leq n$, as $k \geq(n+1) / 2$. If $|c_0| \leq r_{0}$ we are done; otherwise, we get $|c_0|-r_{0} \leq \sqrt{n}$, which finishes the proof.

\medskip

We note that using the left-hand side of \eqref{eqn:3.8} for $r_{0} \leq 1/(2 \sqrt{n})$ we get $|c_0| \geq 1/(2 \sqrt{n})$. This means that in a small enough circle at the origin there are at most $n/2$ of the roots of $G(X, Z)$:

\begin{Cor} There cannot be more than half of the roots of $G(X, Z)$ in a circle of radius $\ds r_0\leq \frac{1}{2\sqrt{n}}$.
\end{Cor}

Now we can return to the binary form $F(X, Z)$ with covariant point $z(F)=(t, u)$. We will use the results above in order to find a bound for $t$, depending on the center $c$ and the radius $r$ of the circle containing more than $n/2$ of the roots.

\begin{Thm}\label{thm:bound_more} If $F(X, Z)$ is a real binary form with more than $n/2$ roots on a disk centered at $c$ and radius $r$, then we have the following bound for the distance $|t-c|$:
\begin{equation}\label{eqn:3.10}
|t-c| \leq(2n+3) r.
\end{equation}
\end{Thm}

\noindent{\bf Proof:} As before, let $\ds\beta_{i}=\frac{t-\alpha_{i}}{u}$ be the roots of $G(X, Z)$, and let $\ds c_0=\frac{t-c}{u}$, $\ds r_{0}=\frac{r}{u}$ be the center and radius of the circle containing the $k>n/2$ roots $\beta_1, \dots, \beta_k$. Then for $1\leq i\leq k$ we have $\ds \frac{|t-\alpha_{i}|}{u} \leq|c_0|+r_{0}$, and further from the right-hand side of \eqref{eqn:3.8} we must have $\ds \frac{|t-\alpha_{i}|}{u} \leq|c_0|+r_{0} \leq \sqrt{n}+2 r_{0}$. Using the bound $u\leq 2r\sqrt{n}$ from Proposition \ref{prop:3.1}, we get
\begin{equation}\label{eqn:3.11}
|t-\alpha_{i}| \leq u\left(\sqrt{n}+2 r_{0}\right) \leq(2n+2)r.
\end{equation}
Finally adding $|\alpha_{i}-c| \leq r$ to the above inequality and applying the triangle inequality we obtain \eqref{eqn:3.10}.

\medskip

\section{Inequalities for precisely $n/2$ roots on a disk of small radius}\label{sec:4}

Let $n\geq 4$ be an even integer. From now on, for the Sections \ref{sec:4}-\ref{sec:6} we assume that we have precisely $n/2$ roots are on a disk of small radius $r_1\leq\epsilon$, centered at $c_1$. The remaining roots will also be contained on a disk of minimum radius $r_2$, centered at $c_2$. Denote the roots in the first circle by $\alpha_1, \dots, \alpha_{n/2}$ and the roots in the second circle by $\alpha_1', \dots, \alpha_{n/2}'$. We will assume throughout the paper that $r_1\leq r_2$.

We will use the notation:
\[
d_i=|t-c_i|, i=1, 2.
\]
Before we continue, we note the easy triangle inequalities that will be used extensively in the proofs:
\begin{equation}\label{eqn:tr1}
|r_1-d_1|\leq |t-\alpha_i|\leq r_1+d_1, \\
\end{equation}
\begin{equation}\label{eqn:tr2}
|r_2-d_2|\leq |t-\alpha'_i|\leq r_2+d_2,  \ \ i=1, 2. 
\end{equation}

In the current section we give several inequalities for $u$ and the quotient $r_1/r_2$ that will be used in the following sections, with the ultimate goal to try to find a conditions for when $u$ is close to zero. In particular, in most cases this will be equivalent to the ratio $r_1/ r_2$ being close to zero. The inequalities we prove in Propositions \ref{prop:4.3}-\ref{prop:4.7} are the main ingredients in showing the Theorems \ref{thm:intro_2}, \ref{thm:intro_3} and \ref{thm:intro_4}(2) from the Introduction.  

We present first some easy upper and lower bounds for $u$, depending on the radii $r_1, r_2$ and on the distances from the coordinate $t$ to the centers $c_1$ and $c_2$.
 
\begin{Prop}\label{prop:4.3} We have $u \leq \sqrt{n}\left( d_1+r_1\right)$ and $u \leq \sqrt{n}\left( d_2+r_2\right)$.
\end{Prop}
\noindent{\bf Proof:} This is an immediate consequence of the Corollary \ref{cor:2.2} (1): as in the circle of center $t$ and radius $ d_1+r_1$ there are $n/2$ roots, the radius must be $\geq u / \sqrt{n}$. Similarly for the second inequality.

\begin{Prop}\label{prop:4.4} The following inequality takes place:
\begin{equation}\label{eqn:4.1}
|\left( d_1-r_1\right)\left( d_2-r_2\right)| \leq u^{2} \leq|\left( d_1+r_1\right)\left( d_2+r_2\right)|.
\end{equation}
\end{Prop}

\noindent{\bf Proof:} Using the triangle inequalities \eqref{eqn:tr1} and \eqref{eqn:tr2}, we get the inequalities:
\[
\begin{aligned}
& \frac{1}{1+\left(d_{1}+r_1\right)^{2} / u^{2}} \leq \frac{1}{1+|t-\alpha_{i}|^{2} / u^{2}} \leq \frac{1}{1+\left(d_{1}-r_1\right)^{2} / u^{2}} \\
& \frac{1}{1+\left(d_{2}+r_2\right)^{2} / u^{2}} \leq \frac{1}{1+|t-\beta_{i}|^{2} / u^{2}} \leq \frac{1}{1+\left(d_{2}-r_2\right)^{2}/u^{2}}
\end{aligned}
\]
We sum up the inequalities for all $1\leq i\leq n/2$, and together with the equation \eqref{eqn:1.1} we get:
\[
\frac{1}{u^{2}+\left(d_{1}+r_1\right)^{2}}+\frac{1}{u^{2}+\left(d_{2}+r_2\right)^{2}} \leq \frac{1}{u^{2}} \leq \frac{1}{u^{2}+\left(d_{1}-r_1\right)^{2}}+\frac{1}{u^{2}+\left(d_{2}-r_2\right)^{2}}.
\]
Easy calculations imply the result.

\medskip

Now we will find bounds for $r_1/r_2$ in the cases of $d_1>r_1$ and $d_2>r_2$. The proof is based on the interpretation of the roots in the hyperbolic plane.

\begin{Prop}\label{prop:4.7}
\begin{enumerate}
	\item For $ d_1>r_1$, we have:
\begin{equation}\label{eqn:4.2}
\frac{3}{n}\frac{\left( d_1-r_1\right)^{2}+u^{2}}{\left( d_2+r_2\right)^{2}+u^{2}} \leq \frac{r_1}{r_2}.
\end{equation}

\item For $ d_2>r_2$, we have:
\begin{equation}\label{eqn:4.5}
\frac{3}{n}\frac{\left( d_2-r_2\right)^{2}+u^{2}}{\left( d_1+r_1\right)^{2}+u^{2}} \leq \frac{r_2}{r_1}. 
\end{equation}

\end{enumerate}
\end{Prop}

\noindent{\bf Proof:} As in section \ref{sec:3.1}, for $1 \leq i \leq n/2$, we denote by $M_i$ (resp. $ M_i'$) the $n/2$ coordinates on the unit sphere that give $\ds\frac{t-\alpha_{i}}{u}$ (resp. $\ds\frac{t-\alpha'_{i}}{u}$) under the stereographic projection, and denote by $\overrightarrow{v_{i}}=\overrightarrow{OM_i}$ (resp. $\vec{w}_{i}$) the $n/2$ vectors on the unit sphere. 

We will prove the first part, with part (2) being proved similarly. As we assume $d_1>r_1$, all $\overrightarrow{v_{i}}$ are in the same half-space. The greatest angle between any two of the vectors $\left(\overrightarrow{v_{i}}, \overrightarrow{v_{j}}\right)$ is attained for some roots, say $\alpha_1$ and $\alpha_2$. Call this angle $\theta$ and using the dot product of $\vec{v}_1\cdot \vec{v}_2$ as in \eqref{eqn:3.4}, we get:

\begin{equation}
\sin \frac{\theta}{2}
=
\frac{|\alpha_{1}-\alpha_{2}|u}{\sqrt{\left(|t-\alpha_{1}|^{2}+u^{2}\right)\left(|t-\alpha_{2}|^{2}+u^{2}\right)}}.\end{equation}

 As $\alpha_1, \alpha_2$ are on the same circle we have $|\alpha_{1}-\alpha_{2}| \leq 2 r_1$, and the triangle inequality gives $|t-\alpha_{i}| \geq\left|d_1-r_1\right|$, for all $i=1, 2$. Using this we find the bound:
\begin{equation}\label{eqn:bound1}
\sin \frac{\theta}{2} \leq \frac{2 r_1 u}{\left( d_1-r_1\right)^{2}+u^{2}}.
\end{equation}

As the roots $\alpha'_{1}, \cdots \alpha'_{n/2}$ can be inscribed in a circle of minimal radius $r_2$, there are at least two points on the boundary of the circle, say $\alpha'_{1}, \alpha'_{2}$, such that $|\alpha'_{1}-\alpha'_{2}| \geq r_2 \sqrt{3}$. Let $\theta'$ be the angle between the corresponding two vectors. We use again the dot product from \refeq{eqn:3.4} and the triangle inequality $|t-\alpha'_{i}| \leq d_2+r_2$ for $i=1, 2$ to obtain:

\begin{equation}\label{eqn:bound2}
\sin \frac{\theta'}{2} 
=
\frac{|\alpha'_{1}-\alpha'_{2}|u}{\sqrt{\left(|t-\alpha'_{1}|^{2}+u^{2}\right)\left(|t-\alpha'_{2}|^{2}+u^{2}\right)}}
\geq 
\frac{r_2 u \sqrt{3}}{\left( d_2+r_2\right)^{2}+u^{2}}
\end{equation}

We note now:

\begin{enumerate}
	\item the minimum possible norm of the vectors $\ds \sum\limits_{i=1}^{n/2}\vec{v_i}$ is $\ds \frac{n}{2} \cos \frac{\gamma}{2}$, where $\gamma$ is the angle of the minimal cone at $O$ containing all the vectors $\vec{v}_1, \dots, \vec{v}_{n/2}$. This is proved in Lemma \ref{lem:3.2}(2).
	\item the maximum possible norm of the vectors $\ds\sum\limits_{i=1}^{n/2}\vec{w_i}$ is $\ds \frac{n}{2}-2+2 \cos \frac{\theta'}{2}$.  This is due to the fact that the sum is maximized when all the vectors $\overrightarrow{v_{3}}, \cdots, \overrightarrow{v_{n/2}}$ have the same direction as the vector $\overrightarrow{v_{1}}+\overrightarrow{v_{2}}$.

\end{enumerate}

Moreover, as the sum of all the unit vectors is $\vec{0}$, the two norms are equal, thus 
\[
\frac{n}{2} \cos \frac{\gamma}{2} \leq \frac{n}{2}-2+2 \cos \frac{\theta'}{2},
\]
which can be rewritten as $\ds (n/2)^{2} \sin^{2} (\gamma/2)-4 \sin^{2}(\theta'/2) \geq 4(n/2-2)(1-\cos(\theta'/2))$. As the right-hand side is clearly non-negative, this implies $\ds \frac{n}{4} \sin(\gamma/2) \geq \sin (\theta'/2)$ and as $\ds\sin(\theta/2) \geq \frac{\sqrt{3}}{2} \sin(\gamma/2)$ from Lemma \ref{lem:3.2}(1), we obtain  
\begin{equation}\label{eqn:angles}
\frac{n/2}{\sqrt{3}} \sin \frac{\theta}{2} \geq \sin \frac{\theta'}{2}.
\end{equation}

Finally the inequalities \eqref{eqn:bound1} and \eqref{eqn:bound2} together with \eqref{eqn:angles} give:
\[
\frac{n/2}{\sqrt{3}} \frac{2 r_1 u}{\left( d_1-r_1\right)^{2}+u^{2}} \geq \frac{n/2}{\sqrt{3}} \sin \frac{\theta}{2} \geq \sin \frac{\theta'}{2} \geq \frac{r_2 u \sqrt{3}}{\left( d_2+r_2\right)^{2}+u^{2}}
\]
which implies the inequality. Case (2) is proved similarly.

\section{Bounds when $u$ is small}\label{sec:5}

In this section we will assume that $u$ is small and using this we will find bounds for the ratio $r_1/r_2$, the product $r_1r_2$ and the distance $d_1=|t-c_1|$.

\medskip

Recall that we have precisely $n/2$ roots $\alpha_1, \dots, \alpha_{n/2}$ on a disk of radius $r_1$, centered at $c_1$, and the other $n/2$ roots $\alpha'_1, \dots, \alpha'_{n/2}$ are on a disk of radius $r_2$, centered at $c_2$. We first note the important observation:

\begin{Prop}\label{prop:4.1} If $u\leq \epsilon$, then either:

	\begin{enumerate}[(i)]
		\item more than $n/2$ roots are on a disk of radius $2\epsilon$
		\item the centers $c_1$ and $c_2$ of the two disks are real numbers, or
		\item $r_1=r_2$ and $c_1=\overline{c_2}$.
		
	\end{enumerate}
\end{Prop}
\noindent{\bf Proof:} We will call $\C_i$ the circle centered at $c_i$ of radius $r_i$ for $i=1,2$. As we are working with real forms, the complex conjugate of a root is also a root of the binary form. We note that if both $\alpha_{i}$ and $\overline{\alpha}_{i}$ are on the circumference of $\C_1$, then its center $c_1$ is real. 

Suppose that $c_1$ is not real. 

\begin{itemize}
	\item If $|\text{Im}(c_1)|>r_1$, then we can reflect the circle along the real axis and obtain a circle centered at $\overline{c_1}$ of radius $r_1$ containing the roots $\overline{\alpha_1}, \dots, \overline{\alpha_{n/2}}$, which are also the roots of the real binary form $F(X, Z)$. Thus we get $r_1=r_2$ and $c_2=\overline{c_1}$, which is case (iii).
	
	\item If $|\text{Im}(c_1)|\leq r_1$, then the circle of radius $2r_1$ centered at $\text{Re}(c_1)$ contains the roots $\alpha_1, \dots, \alpha_{n/2}$ as well as their complex conjugates. If we got more than $n/2$ such roots in the circle of radius $2r_2$, we get case (i). 

Otherwise, we have the same $n/2$ roots on the two circles of radius $r_1$ centered at $c_1$ and $\overline{c_1}$. Thus the roots $\alpha_1, \dots, \alpha_{n/2}$ are inside the intersection of the two circles, and we can find a smaller circle containing all the roots. To see this, denote by $x_1$ and $x_2$ the points of intersection of these two circles. Then the circle centered at $\ds\frac{x_1+x_2}{2}$  of radius $\ds\frac{|x_2-x_1|}{2}<r$ contains all the roots $\alpha_1, \dots, \alpha_{n/2}$, a contradiction. Thus $c_1\in \RR$.
	
\[
\begin{tikzpicture} 
    
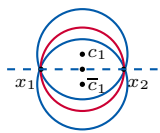
    \draw[cyan!30!blue, thick]   (0,0.2) circle[radius=6mm]  node [yshift=-14mm, xshift=1mm, color=black] {Figure 2: smaller circle containing the roots};

   \fill   (0, 0.2) circle[radius=1pt]  node [above,font=\tiny, yshift=-2mm, xshift=2mm] {$c_1$};
   
    
     \draw[cyan!30!blue, thick]   (0, -0.2) circle[radius=6mm];

      \draw[red!80!blue, thick]   (0,0) circle[radius=5.55mm];
      
      \fill   (-0.55, 0) circle[radius=1pt]  node [above,font=\tiny, xshift=-2mm, yshift=-4mm] {$x_1$};
      
      \fill   (0.55, 0) circle[radius=1pt]  node [above,font=\tiny, xshift=2mm, yshift=-4mm] {$x_2$};
      
      \fill   (0, 0) circle[radius=1pt];

      \fill   (0, -0.2) circle[radius=1pt]  node [above,font=\tiny, yshift=-2.5mm, xshift=2mm] {$\overline{c}_1$};

\draw[cyan!30!blue,dashed, thick]  (-1, -0) -- (1, -0);

    \end{tikzpicture}
\]

A similar argument implies $c_2\in \RR$, as the second disk with circumference $\C_2$ also contains only real roots and pairs of conjugate roots. The smallest disk containing all these roots will have as center a real number.

\end{itemize}

\subsection{Bounds for $r_1/r_2$ when $u$ is small}\label{sec:4.3}
In this section we are showing that when $u$ is small and the two centers $c_1$ and $c_2$ are far away from each other, then the quotient of the two radii $r_1/r_2$ is also close to $0$, i.e. the radius $r_1$ is much smaller than the radius $r_2$.

More precisely, we show:


\begin{Thm}\label{thm:usmall_quot} Let $\ds \epsilon< \frac{1}{n(4 n^{2}+1)}$.
Assuming:
\begin{itemize}
	\item $r_1 \leq \epsilon$ and $u \leq \epsilon$, 
	\item  $|c_1-c_2| \geq \max(4r_2, 4\sqrt{\epsilon})$,
\end{itemize}
then we can bound the ratio $\ds \frac{r_1}{r_2} \leq \frac{10 n \epsilon}{3}.$
\end{Thm}

We will also use these results in the following section in order to show that $t$ and $c_1$ are close to each other.

\smallskip

We will take two cases depending on the size of $r_2$. The first case we consider is of both circles having a small radius, while the two centers being not too close:
\[
\begin{tikzpicture} 
    
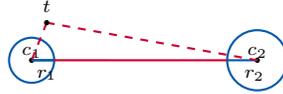
    \draw[cyan!30!blue, thick]   (0,0) circle[radius=3mm] node [yshift=-10mm, xshift=16mm, color=black] {Figure 3: $r_1$ and $r_2$ small, centers far away};
    
   \fill   (0,0) circle[radius=1pt]  node [above,font=\tiny, yshift=-1mm] {$c_1$};
\draw   (0,0) -- node [above, font=\tiny, xshift=-3mm, yshift=-4mm] {$r_1$} (1,0);
\draw[cyan!30!blue, thick] (3,0) circle[radius=4mm] node [text=black,yshift=-2mm, xshift=-0.4mm, font=\tiny] {$r_2$};

 \fill   (3,0) circle[radius=1pt]  node [above,font=\tiny,  yshift=-1mm] {$c_2$};

\draw[cyan!30!blue, thick]   (0,0) -- (0.3 ,0);

\draw[red!80!blue,thick]  (0.3, 0) -- (2.6,0);

\draw[cyan!30!blue, thick]  (2.6, 0) -- (3,0);

\draw[cyan!30!blue, thick]  (0, 0) -- (0.3,0);

\draw[red!80!blue,dashed, thick]  (0.2, 0.5) -- (0,0);
\draw[red!80!blue,dashed, thick]  (0.2, 0.5) -- (3,0);

  \fill   (0.2,0.5) circle[radius=1pt]  node [above,font=\tiny, yshift=0mm] {$t$};

    \end{tikzpicture}
\]

\begin{Prop}\label{prop:4.9}
Let $\ds\epsilon<\frac{1}{n(4 n^{2}+1)}$. 
Assuming $u\leq \epsilon$ and:
\begin{itemize}
	\item  $r_1 \leq \epsilon$, $r_2\leq \sqrt{\epsilon}$,
	\item $|c_1-c_2|>\epsilon+\sqrt{\epsilon}+r_1+r_2$,
	then we can bound the ratio:
\end{itemize} 
\[
\frac{r_1}{r_2} \leq \frac{10 n \epsilon}{3}.
\]
\end{Prop}


\noindent{\bf Proof:} From the left-hand side of the inequality \eqref{eqn:4.1}, as $u \leq \epsilon$, we have:
\[
|d_1-r_1| |d_2-r_2| \leq \epsilon^2,
\]
thus $|d_1-r_1| \leq \epsilon$ or $|d_2-r_2| \leq \epsilon$. Also, note that we cannot have $d_1 \leq r_1$ and $ d_2 \leq r_2$ simultaneously, as the triangle inequality $|c_1-c_2|\leq |t-c_1|+|t-c_2|$ would imply $|c_1-c_2| \leq \epsilon+\sqrt{\epsilon}$.

We will use the inequalities of Proposition \ref{prop:4.7} to get either a bound for $r_1/r_2$, or a contradiction. We first break the case $|d_2-r_2| \leq \epsilon$ into two cases depending on whether $d_1\leq r_1$ or $d_1>r_1$:

\begin{enumerate}
  \item Case $|d_2-r_2| \leq \epsilon$ and $d_1 \leq r_1$. As mentioned above, the condition $d_1 \leq r_1$ implies that we must have $d_2>r_2$, which puts us in the conditions of Proposition \ref{prop:4.7} and \eqref{eqn:4.5} gives us:
  \[
  \frac{r_1}{r_2} 
  \leq \frac{n}{3}\frac{\left(d_1+r_1\right)^{2}+u^{2}}{\left(d_2-r_2\right)^{2}+u^{2}}.
  \]
  Using $d_1+r_1\leq 2r_1\leq 2\epsilon$, $u\leq \eps$, we can bound the numerator by $\left(d_1+r_1\right)^{2}+u^{2}\leq 5\epsilon^2$. For the denominator we use the triangle inequality $|d_2-r_2|\geq d_2-r_2 \geq |c_1-c_2|-d_1-r_2 \geq \sqrt{\eps}$, where at the last step we used the hypothesis $|c_1-c_2|\geq r_1+r_2+\epsilon+\sqrt{\epsilon}$ and $r_1\geq d_1$. Thus we get a bound for the ratio: 
  \[
  \frac{r_1}{r_2} \leq \frac{5 n \epsilon}{3}.
  \]
  
  \item Case $|d_2-r_2| \leq \epsilon$ and $d_1>r_1$. Then we can apply the inequality \ref{eqn:4.2} and we get 
  \begin{equation}\label{eqn:case2}
  \frac{r_2}{r_1} 
  \leq \frac{n}{3}\frac{\left(d_2+r_2\right)^{2}+u^{2}}{\left( d_1-r_1\right)^{2}+u^{2}}.
  \end{equation}
  We bound the numerator using $d_2\leq r_2+\epsilon\leq 2\sqrt{\epsilon}$ and $r_2, u\leq \sqrt{\epsilon}$ to get $\left(d_2+r_2\right)^{2}+u^{2}\leq 5\epsilon$. For the denominator, we use $ d_1-r_1 \geq |c_1-c_2|-d_2-r_1 \geq|c_1-c_2|-r_2-\epsilon-r_1\geq \sqrt{\epsilon}$, in order to get $\ds \frac{r_2}{r_1} \leq \frac{10 n}{3}$. 
  
As $r_1\leq \epsilon$, then we got the lower bound $\ds r_2\leq \frac{10 n\epsilon}{3}$. We will keep lowering the bound for $r_2$ until we get a contradiction such as $r_2<r_1$. Using the new bound for $\ds r_2<\frac{10 n\epsilon}{3}$ in the inequality \ref{eqn:case2}, we get
\[
\frac{r_2}{r_1} \leq \frac{n\left((20 n+3)^{2}+9\right) \epsilon}{27},
\]
which implies $\ds r_2 \leq \frac{n\left((20 n+3)^{2}+9\right) \epsilon^{2}}{27} \leq n \epsilon$, for $\ds \epsilon \leq \frac{27}{(20 n+3)^{2}+9}$. Using the new bound $r_2\leq n\epsilon$ again in inequality \eqref{eqn:case2}, we get
\[
\frac{r_2}{r_1} \leq \frac{n\left((2 n+1)^{2}+1\right) \epsilon}{3}<1,
\]
for $\ds\epsilon<\frac{3}{n\left((2 n+1)^{2}+1\right)}$, contradiction.

\end{enumerate}
Now we break the case of $|d_1-r_1| \leq \epsilon$ into two cases depending on whether $d_2\leq r_2$ or $d_2>r_2$.
\begin{enumerate}
  \item[(3)] Case $|d_1-r_1| \leq \epsilon$ and $d_2 \leq r_2$.  As before, the condition $d_2 \leq r_2$ implies that we must have $ d_1>r_1$ and we can apply the inequality \eqref{eqn:4.2} to get 
  \begin{equation}\label{eqn:case3}
  \frac{r_2}{r_1} \leq \frac{n}{3}\frac{\left(d_2+r_2\right)^{2}+u^{2}}{\left(d_1-r_1\right)^{2}+u^{2}}.
  \end{equation}
  As before, we find the bounds $\left(d_2+r_2\right)^{2}+u^{2}\leq 5\epsilon$ and $\left(d_1-r_1\right)\geq \sqrt{\epsilon}$. This implies $\ds\frac{r_2}{r_1} \leq \frac{5 n}{3}$, and thus $\ds r_2 \leq \frac{5 n \epsilon}{3}$. Applying again \eqref{eqn:case3} for the new bound $\ds r_2 \leq \frac{5 n \epsilon}{3}$, we get
\[
\frac{r_2}{r_1} \leq \frac{n\left(100 n^{2}+9\right) \epsilon}{27}.
\]
Then $\ds r_2 \leq \frac{n\left(100 n^{2}+9\right) \epsilon^{2}}{ 27} \leq n \epsilon$ for $\ds\epsilon \leq \frac{27}{100 n^{2}+9}$. Finally, applying once more \eqref{eqn:case3} with $r_2 \leq n \epsilon$ we obtain 
\[
\frac{r_2}{r_1} \leq \frac{n\left(4 n^{2}+1\right) \epsilon}{3}<1
\] 
for $\ds\epsilon<\frac{1}{n\left(4 n^{2}+1\right)}$, a contradiction.

  \item[(4)] Case of $|d_1-r_1| \leq \epsilon$ and $ d_2>r_2$. We can apply the inequality \eqref{eqn:4.5} and we get 
  \[
  \frac{r_1}{r_2} \leq \frac{n}{3}\frac{\left(d_1+r_1\right)^2+u^{2}}{\left( d_2-r_2\right)^{2}+u^2},
  \]
 Here we bound $\left(d_1+r_1\right)^2+u^{2}\leq \left(\epsilon+2r_1\right)^2+u^2\leq  10\epsilon^2$ and $\ds  d_2-r_2 \geq |c_1-c_2|-d_1-r_2 \geq |c_1-c_2|-r_1-\epsilon-r_2\geq \sqrt{\epsilon}$, and get 
 \[
 \frac{r_1}{r_2} \leq \frac{10 n \epsilon}{3}.
 \]
\end{enumerate}

\medskip

With the help of the above proposition, we can now prove a result also for $r_2>\sqrt{\epsilon}$.

\[
\begin{tikzpicture} 
    
    \draw[cyan!30!blue, thick]   (0,0) circle[radius=3mm] node [yshift=-13mm, xshift=16mm, color=black] {Figure 4: $r_1$ small, centers far away };
   \fill   (0,0) circle[radius=1pt]  node [above,font=\tiny, yshift=-1mm] {$c_1$};
\draw   (0,0) -- node [above, font=\tiny, xshift=-3mm, yshift=-4mm] {$r_1$} (1,0);
\draw[cyan!30!blue, thick] (4,0) circle[radius=10mm] node [text=black,yshift=-2mm, xshift=-6mm, font=\tiny] {$r_2$};

 \fill   (4,0) circle[radius=1pt]  node [above,font=\small] {$c_2$};

\draw[red!80!blue,thick]   (0.3,0) -- (3,0);

\draw[cyan!30!blue, thick]  (0, 0) -- (0.3,0);

\draw[cyan!30!blue, thick]  (3, 0) -- (4,0);

\draw[red!80!blue, thick, dashed]  (0.2, 0.5) -- (4,0);

\draw[red!80!blue, thick, dashed]  (0.2, 0.5) -- (0,0);

  \fill   (0.2,0.5) circle[radius=1pt]  node [above,font=\tiny, yshift=0mm] {$t$};
\
    \end{tikzpicture}
\]

\begin{Prop}\label{prop:4.10} Let $\ds\epsilon<\frac{1}{n(4 n^{2}+1)}$. Assuming $u \leq \epsilon$ and:
\begin{itemize}
	\item $r_1 \leq \epsilon, r_2>\sqrt{\epsilon}$, 
	\item  $|c_1-c_2| \geq 2 r_2+r_2 \sqrt{\epsilon}+r_1$
\end{itemize}
then we can bound the ratio:
\[
\frac{r_1}{r_2} \leq \frac{10 n \epsilon}{3}.
\] 
\end{Prop}

\noindent{\bf Proof:} If $\ds F(X, Z)=a_0\prod_{i=1}^{n/2}(X-\alpha_iZ)(X-\alpha'_iZ)$ is a real binary form with covariant point $z(F)=(t, u)$, then $\ds G(X, Z)=a_0\prod_{i=1}^{n/2}(X-l\alpha_iZ)(X-l\alpha'_iZ)$ is a real binary form with covariant point $z(G)=(lt, lu)$. Clearly, the roots $l\alpha_i$ are on a circle centered at $c'_1=lc_1$ of radius $r'_1=lr_1$, and the roots $l\alpha'_i$ are on a circle centered at $c'_2=lc_2$ of radius $r'_2=lr_2$.

We take $\ds l=\frac{\sqrt{\epsilon}}{r_2}<1$ and then $r_2'=\sqrt{\epsilon}$. Thus we have:
\begin{itemize}
	\item  $r_1'=lr_1 \leq r_1 \leq \epsilon$,
	\item $r_1' \leq r_2'=\sqrt{\epsilon}$ and $u^{\prime} < u \leq \epsilon$.
	\item $\ds |c_1'-c_2'|=\frac{\sqrt{\epsilon}}{r_2}|c_1-c_2|\geq  2 \sqrt{\epsilon}+\epsilon +\frac{r_1\sqrt{\epsilon}}{r_2}=\sqrt{\epsilon}+\epsilon+r_2'+r_1'.
$
\end{itemize}
Thus we are in the conditions of Proposition \ref{prop:4.9} and we obtain
\[
\frac{r_1'}{r_2'}=\frac{r_1}{r_2} \leq \frac{10 n \epsilon}{3}.
\]

\bigskip

Note that in the above proposition we do not need any upper bound for $r_2$. 

\smallskip

The proof of the Theorem \ref{thm:usmall_quot} follows immediately from the last two propositions:

\begin{itemize}
	\item If $r_2\leq \sqrt{\epsilon}$, then $4\sqrt{\epsilon}\geq 4r_2$ and $|c_1-c_2|\geq 4\sqrt{\epsilon}>r_1+r_2+\epsilon+\sqrt{\epsilon}$, thus the result follows from Proposition \ref{prop:4.9}.
	
	\item If $r_2>\sqrt{\epsilon}$, then $4r_2>4\sqrt{\epsilon}$ and $|c_1-c_2|\geq 4r_2>r_2+r_2\sqrt{\epsilon}+r_1$, and the result follows from Proposition \ref{prop:4.10}

\end{itemize}

The remaining case of $|c_1-c_2|\leq \max(4\sqrt{\epsilon}, 4r_2)$ will be treated in  Section \ref{sec:usmall_r1r2}.

\medskip

\subsection{Bounds for $d_1$ when $u$ close to $0$}\label{sec:4.4} In this section, assuming that $u$ is small, we show that $t$ can be approximated by the center $c_1$ of the circle i.e. we bound the value of $d_1=|t-c_1|$ by a value smaller than $1/2$. This result can be used in the algorithm of Stoll and Cremona \cite{SC} using Lemma \ref{lem:u_increase} in Section \ref{sec:algorithm}.


We first present the following lemma, that does not need the assumption of $u$ small. We will use this result below, as well as in Section \ref{sec:6} in order to prove Proposition \ref{prop:4.11}.

\begin{Lem}\label{lem:d2r2_low} Let $\ds\epsilon<\frac{1}{100n(2n+3)^2}$. Assuming that:
	\begin{itemize}
		\item $|c_1-c_2|\geq 2r_2$
		\item $d_1>r_1$
		\item $\ds\frac{r_1}{r_2}<\frac{10n\epsilon}{3}$,
	\end{itemize}
then $|d_2-r_2|\geq r_2/2$ and $d_1\leq u+r_1$.	

\end{Lem}

{\bf Proof:} As $d_1>r_1$, we can apply the inequality \eqref{eqn:4.2} and combined with $\ds r_1/r_2\leq 10n\epsilon/3$ we get
 \[
 \frac{3}{n} \cdot \frac{(d_1-r_1)^{2}+u^{2}}{(d_2+r_2)^{2}+u^{2}}\leq \frac{r_1}{r_2} \leq \frac{10n\epsilon}{3}.
\]

	Multiplying the denominator and rearranging the terms, we have $\ds\frac{10n^2\epsilon}{9} (d_2+r_2)^{2} \geq (d_1-r_1)^{2}+u^2(1-\frac{10n^2\epsilon}{9})\geq (d_1-r_1)^{2}$, thus
\[
\frac{(d_1-r_1)^{2}}{(d_2+r_2)^{2}}\leq \frac{10n^2\epsilon}{9}.
\]
                     
 Taking the square root on both sides, we get $\ds d_1-r_1\leq \frac{n\sqrt{10\epsilon}}{3}(d_2+r_2)$. After dividing by $r_2$, we rewrite:
 \begin{equation}\label{eqn:claim2}
 \frac{d_1}{r_2}-\frac{r_1}{r_2}\leq \frac{n\sqrt{10\epsilon}}{3}\left(\frac{d_2-r_2}{r_2}+2\right) \leq \frac{n\sqrt{10\epsilon}}{3}\left(\left|\frac{d_2-r_2}{r_2}\right|+2\right).
 \end{equation}               
  
Using the triangle inequality twice, we have $d_1=|t-c_1|\geq |c_1-c_2|-|t-c_2|\geq |c_1-c_2|-||t-c_2|-r_2|-r_2$, thus in \ref{eqn:claim2}:
\[
\frac{|c_1-c_2|}{r_2}-\left|\frac{d_2-r_2}{r_2}\right|-1-\frac{r_1}{r_2}\leq \frac{n\sqrt{10\epsilon}}{3}\left(\left|\frac{d_2-r_2}{r_2}\right|+2\right).
\]
After using $r_1/r_2\leq 10n\epsilon/3$ and rearranging the terms, we have:
\[
 \frac{|c_1-c_2|}{r_2}\leq 1+ \frac{2n\sqrt{10\epsilon}}{3}+\frac{10n\epsilon}{3}+\left|\frac{d_2-r_2}{r_2}\right|(1+\frac{n\sqrt{10\epsilon}}{3}).
\]
If $\ds \left|\frac{d_2-r_2}{r_2}\right|<1/2$, then we get $2 r_2 >|c_2-c_1|$ for $\ds\epsilon \leq \frac{1}{160 n^{2}}$. As this is a contradiction, we must have 
\begin{equation}\label{eqn:bound_1/2}
\left|\frac{d_2-r_2}{r_2}\right|\geq 1/2,
\end{equation}
which proves the first result.

Now assume that $d_1-r_1>u$. From the left-hand side of the inequality \eqref{eqn:4.1}, we get:
\[
|d_2-r_2|<u,
\]
thus combined with \eqref{eqn:bound_1/2} gives $u>r_2/2$. This implies $d_2\leq r_2+u\leq 3u$, and further $d_2+r_2\leq 5u$. 

Applying again inequality \ref{eqn:4.2}, we rewrite:
\[
d_2+r_2 \geq u\sqrt{\frac{3r_1}{nr_2}-1}.
\]
Thus we got $\ds 5u \geq d_2+r_2 \geq u\sqrt{\frac{3r_2}{nr_1}-1}$, which is only possible if $\ds \frac{r_1}{r_2}\geq \frac{3}{26n}$. This contradicts the condition $|ds \frac{r_1}{r_2}<\frac{10n\epsilon}{3}$.

\smallskip

Now we are ready to show that $t$ and $c_1$ are close to each other.

\begin{Thm}\label{thm:u_small} Let $\ds\epsilon<\frac{1}{100n(2n+3)^2}$. If $r_1 \leq \epsilon$ and $u\leq \epsilon$, then:
\[
d_1 \leq \frac{1}{2\sqrt{n}}+ \epsilon.
\]
\end{Thm}

{\bf Proof:} (1) Assume first that $r_2>u\sqrt{n}$. Recall from Corollary \ref{cor:2.2}(1) that there are at least $n/2$ roots on a circle centered at $t$ of radius $u\sqrt{n}$. As these cannot be the $n/2$ roots on the circle centered at $c_2$, the circle centered at $t$ of radius $u\sqrt{n}$ contains at least one of the roots $\alpha_i$ of the circle centered at $c_1$. Thus the circle centered at $t$ and radius $u\sqrt{n}$ and the circle centered at $c_1$ and radius $r_1$ intersect and the distance between the centers is smaller than the distance between the radii:
\[
|c_1-t|\leq r_1+u\sqrt{n},
\]
thus $\ds d_1\leq \epsilon+\epsilon\sqrt{n}<\epsilon+\frac{1}{2\sqrt{n}}$.

(2) Assume now that $r_2\leq u\sqrt{n}$, thus $r_2\leq \epsilon \sqrt{n}<\sqrt{\epsilon}$ and the roots $\alpha'_i$ are also on a circle of small radius as well. Note that this implies $r_2\leq \sqrt{\epsilon}$.

(2A) If $|c_1-c_2|<4\sqrt{\epsilon}$, all roots are on a circle of very small radius, and we can put all roots on a circle centered at $c_1$ of radius $4\sqrt{\epsilon}+r_2\leq 5\sqrt{\epsilon}$. Then we are in the conditions of Theorem \ref{thm:bound_more} and we have:
\[
|t-c_1|\leq (2n+3)5\sqrt{\epsilon}\leq \frac{1}{2\sqrt{n}}+\epsilon,
\]
for $\ds\epsilon<\frac{1}{100n(2n+3)^2}$.

(2B) Assume now that $|c_1-c_2|\geq 4\sqrt{\epsilon}\geq \epsilon+\sqrt{\epsilon}+r_2+r_1$. Then we are in the conditions of Proposition \ref{prop:4.9} and we have:
\[
\frac{r_1}{r_2}\leq \frac{10n\epsilon}{3}.
\]

If $d_1<r_1$, we are done. Otherwise we have $d_1>r_1$. As clearly $|c_1-c_2|>2\sqrt{\epsilon}\geq 2r_2$, we are in the conditions of Lemma \ref{lem:d2r2_low} and we have:
\[
d_1\leq r_1+u <r_1+\frac{1}{2\sqrt{n}}.
\]

\medskip

\subsection{Bounds for $r_1r_2$ when $u$ is small}\label{sec:usmall_r1r2}

In this section we treat the remaining case of centers not too far away from each other. 
In this case we do not necessarily expect to have $\ds\frac{r_1}{r_2}<\frac{10n\epsilon}{3}$. Instead, when $u$ is small and the distance of the two centers is not very close to $r_2$, we get a bound for $r_1r_2$, thus not allowing $r_2$ to be too large in comparison to $1/r_1$. 

\[
\begin{tikzpicture} 
    
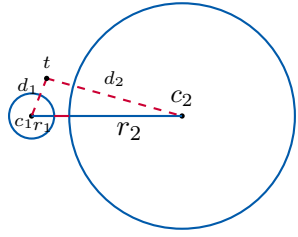
    \draw[cyan!30!blue, thick]   (0,0) circle[radius=3mm] node [yshift=-18mm, xshift=16mm, black] {Figure 5: $r_1$ small; centers not far from each other};
   \fill   (0,0) circle[radius=1pt]  node [above,font=\tiny, yshift=-3mm, xshift=-1mm] {$c_1$};
\draw   (0,0) -- node [above, font=\tiny, xshift=1mm, yshift=-3.5mm] {$r_1$} (0.1,0);
\draw[cyan!30!blue, thick] (2, 0) circle[radius=15mm] node [text=black, yshift=-2mm, xshift=-7mm] {$r_2$};

 \fill   (2, 0) circle[radius=1pt]  node [above,font=\small] {$c_2$};

\draw[cyan!30!blue, thick]  (0, 0) -- (0.3, 0);

\draw[red!80!blue, thick]  (0.3, 0) -- (0.5, 0);

\draw[cyan!30!blue, thick]  (0.5, 0) -- (2, 0);

\draw[red!80!blue, dashed, thick]  (0.2, 0.5) -- (0,0)  node [text=black, yshift=4mm, xshift=-0.5mm, font=\tiny] {$d_1$};;
\draw[red!80!blue, dashed, thick]  (0.2, 0.5) -- (2,0)  node [text=black, yshift=5mm, xshift=-9mm, font=\tiny] {$d_2$};;

  \fill   (0.2,0.5) circle[radius=1pt]  node [above,font=\tiny, yshift=0mm] {$t$};

    \end{tikzpicture}
\]

More precisely, we show:

\begin{Prop}\label{prop:usmall_prod} Let $\ds\epsilon<\frac{1}{100n(2n+3)^2}$ and $u\leq \epsilon$. We are considering the case of:
	\begin{itemize}
		\item $r_1\leq \epsilon$
		
		\item $|c_1-c_2|\leq \max(r_1+r_2+\epsilon+\sqrt{\epsilon}, 2r_2+r_2\sqrt{\epsilon}+r_1)$.
	\end{itemize}
	
Assume the extra hypothesis that $\ds |d_1-r_1|\geq \frac{u}{\sqrt{n}}$. If $|c_1-c_2|$ is not close to the value of $r_2$, meaning that
	\[
	\left|\frac{|c_1-c_2|}{r_2}-1\right| \geq 64n\sqrt{n}\epsilon
	 \]
then we have the bound\[
		r_1r_2\leq \frac{1}{64n^2}.
		\]
\end{Prop}

{\bf Proof:}  Assume first that $r_2\leq \sqrt{\epsilon}$. Then $\ds r_1r_2\leq \epsilon\sqrt{\epsilon}<\frac{1}{64n^2}$, and we are done.

Now assume that $r_2>\sqrt{\epsilon}$. As $\ds |d_1-r_1|\geq \frac{u}{\sqrt{n}}$, we have from the left-half side of the inequality from Proposition \ref{prop:4.4} that $|d_2-r_2| \leq u\sqrt{n}$, which implies 
 \begin{equation}\label{eqn:ineq_d2}
 d_2\leq r_2+u\sqrt{n}.
 \end{equation}
We also note that as $u\leq \epsilon$, from Theorem \ref{thm:u_small} we have 
\begin{equation}\label{eqn:ineq_d1}
d_1\leq r_1+\frac{1}{2\sqrt{n}}.
\end{equation}

We will use these two inequalities to find a bound for $r_1r_2$. We consider the two cases:

(1) $|c_1-c_2|= k r_2$, for $k\geq 1+64n\sqrt{n}\epsilon$.
From the triangle inequality and the inequalities above we have:
\[
kr_2= |c_1-c_2|\leq d_2+d_1 \leq r_2+u\sqrt{n}+\frac{1}{2\sqrt{n}}+r_1 .
\]
We can bound the right-hand side by $r_2+ \frac{1}{\sqrt{n}}$ and get $\ds r_2\leq \frac{1}{(k-1)\sqrt{n}}$. Then 
\[
r_1r_2\leq \frac{\epsilon}{(k-1)\sqrt{n}} \leq \frac{1}{64n^2}
\]
 for $k \geq 1+64n\sqrt{n}\epsilon$.

(2) $|c_1-c_2|= kr_2$ for $k\leq 1-64n\sqrt{n}\epsilon$. From the triangle inequality and the inequalities above we have:
\[
d_2\leq |c_1-c_2|+d_1\leq kr_2+\frac{1}{2\sqrt{n}}+r_1,
\]
As $r_2\leq d_2+u\sqrt{n}$, then 
\[
r_2\leq kr_2+u\sqrt{n}+\frac{1}{2\sqrt{n}}+r_1\leq kr_2+\frac{1}{\sqrt{n}}
\]
thus $r_2\leq \frac{1}{(1-k)\sqrt{n}}$ and
\[
r_1r_2 \leq \frac{\epsilon}{(1-k)\sqrt{n}} <\frac{1}{64n^2}
\]
for $64n\sqrt{n}\epsilon \leq 1-k$.

\section{Showing that $u$ is close to $0$}\label{sec:6}

Recall that we have $n/2$ roots $\alpha_1, \dots, \alpha_{n/2}$ on a disk of small radius $r_1<\epsilon$, centered at $c_1$ and the other $n/2$ roots $\alpha_1', \dots, \alpha'_{n/2}$ on a disk of minimal radius $r_2$, centered at $c_2$. The goal of the section is to show that if either the ratio $r_1/r_2$ is small or the product $r_1r_2$ is bounded by a small enough value, then $u$ must also be small. Combined with the results from the previous section, this will give a necessary and sufficient condition for $u$ to be small. Moreover we give sufficient conditions for $t$ and $c_1$ to be close to each other. 




\subsection{Bounds independent of the size of $u$}\label{sec:bounds_indep}

In this section we will show that when $u$ is not necessarily close to $0$, but the ratio of the two roots $r_1/r_2$ is close to $0$ or the product $r_1r_2$ is relatively small, then we can also approximate the coordinate $t$ by $c_1$.  We show:

\begin{Thm}\label{thm:u_large}
Let $\ds \epsilon \leq\frac{1}{100n(2n+3)^2}$. If $r_1\leq \epsilon$ and either of the conditions:
	\begin{itemize}
		\item $|c_1-c_2|\geq 2r_2$ and $\ds\frac{r_1}{r_2}\leq \frac{10n\epsilon}{3}$, or
		\item  $|c_1-c_2|\leq 2r_2$, $c_1, c_2\in \RR$ and $\ds r_1r_2\leq \frac{3}{64n^2}$,
	\end{itemize}
then we have the bound 
\[
d_1 \leq \frac{1}{2 \sqrt{n}}+r_1.
\]
\end{Thm}
We will show in Section \ref{sec:u_small} that actually in these cases $u$ must be close to $0$. 

\bigskip
We first consider the case of the centers $c_1, c_2$ far away from each other.
\[
\begin{tikzpicture} 
    
    \draw[cyan!30!blue, thick]   (0,0) circle[radius=3mm] node [yshift=-17mm, xshift=16mm, black] {Figure 6: $r_1$ small, centers far away };
   \fill   (0,0) circle[radius=1pt]  node [above,font=\tiny, yshift=-1mm] {$c_1$};
\draw   (0,0) -- node [above, font=\tiny, xshift=-3mm, yshift=-4mm] {$r_1$} (1,0);
\draw[cyan!30!blue, thick] (5,0) circle[radius=15mm] node [text=black,yshift=2mm, xshift=-6mm] {$r_2$};

 \fill   (5,0) circle[radius=1pt]  node [above,font=\small] {$c_2$};

\draw[red!80!blue, thick]  (0.3, 0) -- (3.5,0);

\draw[cyan!30!blue, thick]  (0, 0) -- (0.3,0);

\draw[cyan!30!blue, thick]  (3.5, 0) -- (5,0);

\draw[red!80!blue, dashed, thick]  (0.2, 0.5) -- (0,0);
\draw[red!80!blue, dashed, thick]  (0.2, 0.5) -- (5,0);

  \fill   (0.2,0.5) circle[radius=1pt]  node [above,font=\tiny, yshift=0mm] {$t$};

    \end{tikzpicture}
\]

\begin{Prop}\label{prop:4.11} 
Let $\ds \epsilon \leq \frac{1}{100n(2n+3)^2}$. Assuming the conditions:
\begin{itemize}
	\item $|c_1-c_2| \geq 2 r_2$,
	\item $\ds r_1 \leq \epsilon$ and $\ds \frac{r_1}{r_2} \leq \frac{10 n \epsilon}{3}$, 
\end{itemize} then we have the bound for $d_1$: 
\[
 d_1 \leq \frac{1}{2 \sqrt{n}}+r_1.
\]
\end{Prop}


\noindent{\bf Proof:} Suppose $\ds d_1 \geq \frac{1}{2 \sqrt{n}}+r_1$. Using the hypothesis $r_1/r_2\leq 10n\epsilon/3$ and $|c_1-c_2|\geq 2r_2$, we are in the conditions of Lemma \ref{lem:d2r2_low} and we have:
\begin{equation}\label{eqn:boundr2_2}
|d_2-r_2|\geq r_2/2.
\end{equation}
We will first show that this implies $u\geq 1/n^2$. As from the assumption we have $\ds d_1-r_1\geq  \frac{1}{2\sqrt{n}}$, from the left-hand side of the inequality \eqref{eqn:4.1}:
\[
(d_1-r_1) |d_2-r_2| \leq u^2,
\]
we must have $\ds |d_2-r_2| \leq 2 \sqrt{n} u^{2}$. Together with $|d_2-r_2|\geq r_2/2$ above, this gives us $r_2\leq 4 \sqrt{n} u^{2}$ and thus:
\begin{equation}\label{eqn:lowbound_u2}
d_2+r_2\leq 10u^2\sqrt{n}.
\end{equation}
Applying the inequality \eqref{eqn:4.2}, we get 
\[
d_2+r_2\geq u\sqrt{\frac{3r_1}{n r_2}-1}.
\]
Together with \eqref{eqn:lowbound_u2}, we get:
\[
100 nu^{2} \geq \frac{9}{10 n^{2} \epsilon}-1
\]
and for $\ds\epsilon \leq \frac{9}{10 n(n+100)}$ we get 
\begin{equation}\label{eqn:lowbound_uu}
u \geq \frac{1}{n^{2}}.
\end{equation}
We will go through the same steps once more using the inequality \ref{eqn:lowbound_uu} to get a contradiction. From Proposition \ref{prop:4.3}, we have $\ds  d_1+r_1 \geq \frac{u}{\sqrt{n}}$. As $\ds u \geq \frac{1}{n^{2}}$ and $r_1 \leq \epsilon$, then we can write 
\[
 d_1-r_1 \geq \frac{u}{\sqrt{2 n}}
\]
for $\ds\epsilon \leq \frac{1}{6 n^{2} \sqrt{2 n}}$. Applying again the left-hand side of the inequality \ref{eqn:4.1} we must have 
\[
|d_2-r_2| \leq u \sqrt{2n},
\] 
which together with \eqref{eqn:boundr2_2} gives us $\ds r_2 \leq 2 u \sqrt{2 n}$ and $d_2+r_2 \leq 5 u \sqrt{2 n}$.

Let us apply once more \eqref{eqn:4.2} with the bounds $d_2+r_2 \leq 5 u \sqrt{2 n}$ and $\ds d_1-r_1\geq \frac{u}{\sqrt{2n}}$:
 \[
 \frac{10 n^2 \epsilon}{9} \geq \frac{\left( d_1-r_1\right)^{2}+u^{2}}{\left( d_2+r_2\right)^{2}+u^{2}} \geq \frac{(u / \sqrt{2 n})^{2}+u^{2}}{(5 u \sqrt{2 n})^{2}+u^{2}}=\frac{3(1+1 / 2 n)}{n(50 n+1)},
 \]
 implying $\ds \epsilon \geq \frac{9(1+1 / 2 n)}{10 n^{2}(50 n+1)}$, which is a contradiction.

\bigskip

Now we consider the case of the two centers being possibly closer to each other (see Figure 5), meaning that $|c_1-c_2|\leq 2r_2$. In this case we assume that the product of the two radii is not too large.

Recall from Proposition \ref{prop:4.1} that when $u\leq \epsilon$ is small, we must have one of the conditions: 
\begin{enumerate}[(i)]
	\item more than half of the roots are on a circle of radius $2\epsilon$, centered at $\text{Re}(c_1)$  
	\item $c_1, c_2\in \RR$,
	\item $c_1=\overline{c_2}$ and $r_1=r_2$.
\end{enumerate} 

Assuming $|c_1-c_2|\leq 2r_2$, in case (iii) we have $|c_1-c_2|\leq 2r_1\leq 2\epsilon$, thus all roots are in a circle of radius $2\epsilon$ centered at $\text{Re}(c_1)$. Then in both cases (i) and (iii), from Theorem \ref{thm:bound_more} we get:
\[
|t-\text{Re}(c_1)|\leq 2(2n+3)\epsilon.
\]

Thus the only case left when we could have $u$ small and $|c_1-c_2|\leq 2r_2$ is that of $c_1, c_2\in \RR$. We show below that also in this case $t$ and $c_1$ are close to each other when $r_1r_2$ is bounded by a small enough value.

\begin{Prop}\label{prop:4.13}
Let $\ds\epsilon \leq \frac{3}{104 n^{2}(n+1)}$. Assuming:

\begin{itemize}
	\item $r_1\leq \epsilon$ and $\ds r_1r_2\leq \frac{3}{64 n^{2}}$
	\item $c_1, c_2\in \RR$ and $|c_1-c_2| \leq 2 r_2$,
\end{itemize}	
then we have the bound $\ds d_1 \leq \frac{1}{2 \sqrt{n}}+r_1$.
\end{Prop}

\noindent{\bf Proof:} Suppose $ \ds d_1>r_1+\frac{1}{2 \sqrt{n}}$. Then we can apply the inequality \eqref{eqn:4.2} which we write in the form: 
\begin{equation}\label{eqn:ineq_r1r2}
\frac{3}{n} \cdot \frac{r_2}{r_1}
\leq
\frac{\left( d_2+r_2\right)^{2}+u^{2}}{\left(d_1-r_1\right)^{2}+u^{2}} \leq  4n (\left( d_2+r_2\right)^{2}+u^2),
\end{equation}
where in the last inequality we used the bounds $\ds d_1-r_1 \geq \frac{1}{2\sqrt{n}}$ and $u^2\geq 0$. Our goal is to find bounds depending on $r_2$ on the right hand side, which will lead to a contradiction on the assumption $\ds r_1r_2 \leq \frac{3}{64n^2}$.

First, from Proposition \ref{prop:4.3} we have the inequality $\ds u \leq\left(r_2+ d_2\right) \sqrt{n}$, thus it is enough to bound $d_2+r_2$. In order to do this, note that from \cite{J}, $t$ has to be in the convex hull (in the hyperbolic space) of the roots, thus 
\[
\min \left(c_1-r_1, c_2-r_2\right) \leq t \leq \max \left(c_1+r_1, c_2+r_2\right).
\] 
This implies $d_2=|t-c_2|\leq \max (r_2, |c_1-c_2|+r_1)$, thus $d_2\leq 2r_2+r_1\leq 3r_2$ and
\begin{equation}\label{eqn:r2d2}
 d_2+r_2 \leq 4r_2.
\end{equation}

Then using $u\leq (d_2+r_2)\sqrt{n}$ and \ref{eqn:r2d2} in the inequality \ref{eqn:ineq_r1r2}, we get
\[
\frac{r_2}{r_1} \leq\frac{4n^2}{3}
(d_2+r_2)^2(n+1)\leq 
\frac{64n^2}{3}(n+1)r_2^2,
\]
thus $\ds r_1 r_2\geq \frac{3}{64 n^{2}(n+1)}  $. 
\smallskip

This is not a contradiction yet. However, 
as $r_1\geq\epsilon$, it implies $r_2 \geq 1$ for $\ds \epsilon \leq \frac{3}{64 n^{2}(n+1)}$. Then we can bound $\ds\frac{n}{3}\frac{r_1}{r_2}\leq \frac{n\epsilon}{3}<1$ and looking again at \eqref{eqn:4.2}, we can rewrite 
 \[
 \frac{n}{3} \cdot \frac{r_1}{r_2} \cdot\left( d_2+r_2\right)^{2} \geq u^{2}\left(1-\frac{n}{3} \cdot \frac{r_1}{r_2}\right)+( d_1- \left.r_1\right)^{2} \geq\left( d_1-r_1\right)^{2}.
 \] 
 Thus, using \ref{eqn:r2d2} and the assumption that $d_1-r_1\geq\frac{1}{2\sqrt{n}}$, we get
\[
\frac{n}{3} \cdot \frac{r_1}{r_2} \cdot 16 r_2^{2} \geq \frac{1}{4 n},
\] 
implying $\ds r_1 r_2 \geq \frac{3}{64 n^{2}}$, a contradiction.

\smallskip

\medskip

From Proposition \ref{prop:4.11} and \ref{prop:4.13}, we get Theorem \ref{thm:u_large}.

\subsection{Showing that $u$ is small}\label{sec:u_small}

In this section we will show that under the conditions:
	\begin{itemize}
		\item $\ds \frac{r_1}{r_2}<\frac{10n\epsilon}{3}$, or
		\item $\ds r_1r_2<\epsilon^2$
	\end{itemize}
we obtain that $u$ is small. These are converse statements for Theorem \ref{thm:usmall_quot} and Proposition \ref{prop:usmall_prod}. More precisely, we will show:
\begin{Thm} Let $\ds\epsilon<\frac{1}{100n(2n+3)^2}$ and  $r_1\leq \epsilon$. 

(1)Assume either: 
\begin{itemize}
	\item $r_2\leq \sqrt{\epsilon}$, or 
	\item $r_2> \sqrt{\epsilon}$ and $|c_1-c_2|\geq 2r_2$. 
\end{itemize}	

Then if $\ds \frac{r_1}{r_2}\leq \frac{10n\epsilon}{3}$, we get 
\[
\ds u\leq \frac{20n}{7}\epsilon.
\]

(2) In the remaining case:
\begin{itemize}
	\item  $r_2> \sqrt{\epsilon}$ and $|c_1-c_2|\leq 2r_2$,
\end{itemize}
then if $r_1r_2\leq \epsilon^2$, we get 
\[
u\leq 2\sqrt{\epsilon}.
\]
	
\end{Thm}

We first consider the case of $r_1, r_2$ small (see Figure 3), with the quotient $r_1/r_2$ small. 

\begin{Prop}\label{prop:r2small_u} Assume $r_1\leq \epsilon$ and $r_2\leq \sqrt{\epsilon}$. If $\ds \frac{r_1}{r_2}<\frac{10n\epsilon}{3}$,  then $\ds u<\frac{20n}{7}\epsilon$.

\end{Prop}

{\bf Proof:} We consider the two cases $d_1\leq r_1$ and $d_1>r_1$:

\smallskip
(1) If $d_1\leq r_1$, recall that from Proposition \ref{prop:4.3} we have $u/\sqrt{n} \leq d_1+r_1$. As from the assumption we have $d_1\leq r_1\leq \epsilon $, this gives the bound
\[
u/\sqrt{n}\leq 2r_1\leq 2\epsilon.
\]

(2) If $d_1>r_1$, we will find upper and lower bounds for $d_2-r_2$ with respect to $u$. As $d_1>r_1$, we are in the conditions of Proposition \ref{prop:4.7}  (1) and we rewrite:
\[
\frac{n}{3}\frac{r_1}{r_2}(d_2+r_2)^2 \geq u^2(1-\frac{n}{3}\frac{r_1}{r_2}).
\] 
Thus using the hypothesis $\ds \frac{r_1}{r_2}\leq \frac{10n\epsilon}{3}$, we get $\ds d_2+r_2\geq u \sqrt{\frac{9}{10n^2\epsilon}-1}$, which together with $r_2\leq \sqrt{\epsilon}$ further implies
\begin{equation}\label{eqn:diff_low}
d_2-r_2\geq u \sqrt{\frac{9}{10n^2\epsilon}-1}-2\sqrt{\epsilon}.
\end{equation}
Moreover, from Proposition \ref{prop:4.3} we have $d_1+r_1\geq u/\sqrt{n}$, which implies $d_1-r_1\geq u/\sqrt{n}-2\epsilon$. Combining this with the the left-hand side of the inequality \ref{eqn:4.1} we get the upper bound:
\begin{equation}\label{eqn:diff_upp}
d_2-r_2 \leq \frac{u^2\sqrt{n}}{u-2\epsilon\sqrt{n}}.
\end{equation}

Then \eqref{eqn:diff_low} and \eqref{eqn:diff_upp} imply
\[
\frac{u^2\sqrt{n}}{u-2\epsilon\sqrt{n}}\geq d_2-r_2\geq u \sqrt{\frac{9}{10n^2\epsilon}-1}-2\sqrt{\epsilon},
\]
which simplifies to a bound $\ds u<\frac{20n \epsilon}{7}$.

\medskip

Now we consider $r_2$ not too small and the centers far away (see Figure 4), with the quotient $r_1/r_2$ small.

\begin{Prop}\label{prop:r2notsmall_u} Let $\ds\epsilon<\frac{1}{100n(2n+3)^2}$ and $r_1\leq \epsilon$. Assume the conditions:
\begin{itemize}
	\item  $r_2>\sqrt{\epsilon}$ 
	\item  $|c_1-c_2|\geq 2r_2$.
\end{itemize}

If $\ds\frac{r_1}{r_2}\leq \frac{10n\epsilon}{3}$, then $u<8\sqrt{n}\epsilon$.

\end{Prop}

{\bf Proof:} We split the proof into two cases:

(1) If $d_1\leq r_1$, then from Proposition \ref{prop:4.3}  we have $u/\sqrt{n}\leq d_1+r_1\leq 2\epsilon$.

(2) If $d_1>r_1$, we split the solution into two further cases:

(2A) Case of $d_2\leq r_2$. Using the hypothesis $2r_2\leq |c_1-c_2|$ and the triangle inequality, we have:
\[
2r_2\leq |c_1-c_2| \leq d_2+d_1\leq r_2+d_1,
\]
thus $d_1\geq r_2$. 

Applying Proposition \ref{prop:4.7}(1) and using the bound $r_1/r_2\leq 10n\epsilon/3$, we get:
\[
\frac{10n^2 \epsilon}{9}(d_2+r_2)^2\geq u^2(1-\frac{10n^2 \epsilon}{9})+(d_1-r_1)^2 \geq (d_1-r_1)^2.
\]
Extracting the square root, this implies:
\begin{equation}\label{eqn:bound_diff1}
\frac{n \sqrt{10\epsilon}}{3}(d_2+r_2) \geq d_1-r_1.
\end{equation}
Using the hypothesis of (2A) that $d_2\leq r_2$ and the inequality $d_1\geq r_2$ that we proved above, we get in \ref{eqn:bound_diff1}:
\[
r_2\frac{2n \sqrt{10\epsilon}}{3} \geq r_2-r_1,
\]
thus $\ds r_1 \geq (1-\frac{2n \sqrt{10\epsilon}}{3})r_2\geq \frac{1}{2} r_2$, implying $\ds r_2<\sqrt{\epsilon}$, a contradiction.

(2B) Case of $d_2>r_2$. The goal is to find a lower and an upper bound for $d_2+r_2$ depending only on $u$ and constants. This will allow us to find a bound for $u$ as well. 

 We start by applying Proposition \ref{prop:4.7}(1), which we can rewrite as:
\[
(d_2+r_2)^2\geq u^2\left(\frac{3}{n}\frac{r_2}{r_1}-1\right)+(d_1-r_1)^2 \frac{3}{n}\frac{r_2}{r_1},
\]
and using $(d_1-r_1)^2\geq 0$, we get the lower bound 
\begin{equation}\label{eqn:bound_upd2r2}
d_2+r_2\geq u\sqrt{\frac{2r_2}{nr_1}}.
\end{equation}

Assume that $u>8\epsilon\sqrt{n}$. From Proposition 
\ref{prop:4.3} we have $\ds d_1+r_1 \geq \frac{u}{\sqrt{n}}$, which implies $\ds d_1-r_1 \geq \frac{u}{\sqrt{2n}}$ as $u>8\epsilon\sqrt{n}$. Then from the left-hand side of the inequality of Proposition \ref{prop:4.4}, we get
\begin{equation}\label{eqn:lowbound_u}
d_2-r_2 \leq u\sqrt{2n}.
\end{equation}
We note that we are in the conditions of Lemma \ref{lem:d2r2_low}, thus:
\[
d_2-r_2\geq\frac{r_2}{2}.
\]
Together with \eqref{eqn:lowbound_u}, this implies $r_2\leq 2u\sqrt{2n}$ and $d_2+r_2\leq 5u\sqrt{2n}$. Together with \eqref{eqn:bound_upd2r2} we get:
\[
u\sqrt{\frac{2r_2}{nr_1}}\leq d_2+r_2 \leq 5u\sqrt{2n}, 
\]
which implies $\ds\frac{r_1}{r_2}\geq \frac{1}{25n}$, which contradicts the bound $\ds\frac{r_1}{r_2}\leq \frac{10n\epsilon}{3}$. Thus we must have $u\leq 8\epsilon\sqrt{n}$.

\medskip

Finally, we are left with the case of $r_2$ not too small and the centers close to each other (see Figure 5). From Proposition \ref{prop:usmall_prod} (and Proposition \ref{prop:4.13}), we expect that in order for $u$ to be small we need to have $r_1r_2$ bounded by a small enough value. Taking the product $r_1r_2$ small, we show that $u$ is also small.

\begin{Prop}\label{prop:r1r2_small} Let $\ds\epsilon<\frac{1}{100n(2n+3)^2}$. Assuming the conditions:

	\begin{itemize}
		\item $r_1\leq \epsilon$ and $r_2>\sqrt{\epsilon}$
		\item $|c_1-c_2|<2r_2$
		\item $\ds r_1r_2\leq \frac{3}{64n^2}$,
	\end{itemize}
we have the bound for $u$:
\[
u\leq 4\sqrt{nr_1r_2}.
\]	
If we further have $r_1r_2\leq \epsilon^2$, then $u\leq 4\epsilon\sqrt{n}$.
\end{Prop}

{\bf Proof:} We note first that from the bounds on $r_1$ and $r_2$ we have:
\[
\frac{r_1}{r_2}\geq \frac{\epsilon}{\sqrt{\epsilon}} \geq n^2\epsilon.
\]
Moreover, as $\ds r_1r_2\leq \frac{3}{64n^2}$, Proposition \ref{prop:4.13} implies that $\ds |d_1-r_1|\leq \frac{1}{2\sqrt{n}}$. We take the two cases:

(1) $d_1\leq r_1$. In this case, as from Proposition \ref{prop:4.3} we have $\ds d_1+r_1\geq \frac{u}{\sqrt{n}}$, we get $u\leq 2r_1/\sqrt{n}\leq 2\epsilon/\sqrt{n}$.

(2) $d_1>r_1$.  In this case we can apply the inequality of Proposition \ref{prop:4.7}(1) in order to get:
\begin{equation}\label{eqn:ineq_d2r2u}
d_2+r_2\geq u\sqrt{\frac{3r_2}{nr_1}-1}\geq u \sqrt{\frac{r_2}{nr_1}}.
\end{equation}

We want to find a bound of $d_2$ depending on $r_2$. We apply the triangle inequality, followed by the hypothesis bound of $|c_1-c_2|\leq 2r_2$ and the inequality $d_1 \leq \frac{1}{2\sqrt{n}}+r_1$ of Proposition \ref{prop:4.13} to get:
\[
d_2\leq |c_1-c_2|+d_1\leq 2r_2+(\frac{1}{2\sqrt{n}}+r_1).
\]
We note that $\ds\frac{1}{2\sqrt{n}}+r_1\leq \frac{1}{\sqrt{n}}\leq r_2$, thus $\ds d_2\leq 3r_2$ and together with \eqref{eqn:ineq_d2r2u} we get:
\[
 u \sqrt{\frac{r_2}{nr_1}}\leq d_2+r_2 \leq 4r_2.
\]
Thus $u\leq 4 \sqrt{nr_1r_2}$. If $r_1r_2\leq \epsilon^2$, we immediately get $u\leq 4\epsilon\sqrt{n}$.


\medskip

\begin{rmk}We further remark that a similar argument can prove the more general statement:
\[
d_1-r_1\leq \max\left(\sqrt{r_1r_2\frac{16n}{3}}, \frac{16n(n+1)\epsilon}{3}\right),
\]
when $r_1\leq \epsilon$, $c_1, c_2\in \RR$ and $|c_1-c_2|\leq 2r_2$.
\end{rmk}

\begin{rmk}\label{rmk:bound} We note that in the case 
\begin{itemize}
	\item $r_1\leq \epsilon$, $r_2\leq \sqrt{\epsilon}$
	\item $|c_1-c_2|\leq 2r_2$,
\end{itemize}
we can bound $u$ without any extra condition. As all the roots are on a circle of radius $2\sqrt{\epsilon}$, from Proposition \ref{prop:3.1} we have the bound:
\begin{equation}
u\leq 4\sqrt{n\epsilon}.
\end{equation}

\end{rmk}

\section{The algorithm}\label{sec:algorithm}

Let $F(X, Z)$ be a real binary form with roots $\alpha_1, \dots, \alpha_n$. We assume that the roots are known with sufficient precision. 

We note that one can easily check if at least half of a roots are in the a circle of radius at most $\ds\epsilon$. We first can compute the differences between the roots
\[
|\alpha_j-\alpha_k|\leq 2\epsilon, j \neq k,
\]
for all $j\neq k$ in order to check if they could be on a circle of radius $2\epsilon$.

If we have at least $n/2$ such roots, we find the circle with the smallest radius containing the roots. The best algorithm, due to Nimrod Megiddo (1985) has $O(n)$ complexity. In this way we find  the center $c_1$ and radius $r_1$. We will denote these roots $\alpha_{1}, \ldots, \alpha_{k}$, for $k \geq n/2$. As the distance between any of the two roots is $\leq 2 \epsilon$, we have $r_1 \leq 2 \epsilon/\sqrt{3}$.

The algorithm of Stoll and Cremona \cite{SC} computes an $\SL_2(\ZZ)$-equivalent binary form, with the goal of moving the covariant point in the fundamental domain $\F$ of $\HH$ under the action of $\SL_2(\ZZ)$, as described in the introduction. If $z(F)=(t, u)\in \HH$ with $u$ small, let $m=\lfloor t+1/2\rfloor$ the integer closest to the coordinate $t$,  and apply the transformation:
\[
F'(X, Z):=F(Z, -X-m Z).
\]
The covariant number of the new binary form is 
\[
z(F'):=\left(\frac{-(t-m)}{(t-m)^{2}+u^{2}}, \frac{u}{(t-m)^{2}+u^{2}}\right).
\] 
The purpose of this transformation is to significantly increase the value of $u$ after each step. 

When $z(F)$ is close to the real axis, the number cannot be computed with reasonable accuracy with the usual algorithm. Our goal is to use $\text{Re}(c)$ instead of $t$ in the algorithm, i.e. $m=\lfloor \text{Re}(c)+1/2\rfloor$ instead of $m=\lfloor t+1/2\rfloor$.

From Corollary \ref{cor:usmall}, if $u$ is small there are at least $n/2$ roots in a circle of small radius, say centered at $c$ and radius $r$. We showed in Theorem \ref{thm:bound_more} and Theorem \ref{thm:u_small} that $t$ and $c$ are very close to each other: $|t-c|\leq (2n+3)$ if $k>n/2$ and $|t-c|\leq \frac{1}{2\sqrt{n}}+r$ if $k=n/2$. Note that the same inequalities hold also for $\text{Re}(c)$ instead of $c$.

If we take $m=\lfloor \text{Re}(c)+1/2\rfloor$ in the algorithm, we show below that the value of $u$ will increase as well by a factor of at least $8/7$.

\begin{Lem}\label{lem:u_increase} Let $F(X, Z)$ be a binary form with precisely $k$ roots on a disk of radius $r \leq \epsilon$ and center $c$ and let $m=\lfloor \text{Re}(c)+1/2 \rfloor$.
\begin{enumerate}
	\item If $k>n/2$ and $\ds\epsilon\leq \frac{1}{32(n+1)}$, then: 
\[
\frac{u}{(t-m)^{2}+u^{2}}>2 u
\]
	\item If $k=n/2\geq 2$, $\ds\epsilon \leq \frac{1}{4n^2}$ and $d_1 \leq \frac{1}{2 \sqrt{n}}+\epsilon$ then:
	\[
	\frac{u}{(t-m)^{2}+u^{2}} \geq \frac{8}{7} u.
	\]
\end{enumerate}	
\end{Lem}

\noindent{\bf Proof:} (1) For $k>n/2$, from Proposition \ref{prop:3.1} we have the bound $u \leq 2r \sqrt{n}$.

Also, using the inequality \eqref{eqn:3.10} that $|t-c|\leq (2n+3)r$ and the fact that $m$ is the integer closest to $\text{Re}(c)$, we obtain $|t-m| \leq$ $\ds\frac{1}{2}+|t-\text{Re}(c)| \leq \frac{1}{2}+(2 n+3) r$. The two inequalities imply:
\[
(t-m)^{2}+u^{2} \leq \frac{1}{4}+(2 n+3)^{2} r^{2}+(2 n+3) r+4nr.
\]
As $(2n+3)r<1$, we can bound $(2 n+3)^{2} r^{2}+(2n+3) r+4 n r\leq (8 n+6) r\leq 1/4$. Thus we get $(t-m)^{2}+u^{2}<1/2$.

(2) In the second case, as $\ds |t-m| \leq \frac{1}{2}+ d_1$ and $u \leq \sqrt{n}\left( d_1+r_1\right)$ from Proposition \ref{prop:4.3}, using the hypothesis $d_1\leq\frac{1}{2\sqrt{n}}+\epsilon$ we have 
\[
(t-m)^{2}+u^{2} \leq\left(\frac{1}{2}+\frac{1}{2 \sqrt{n}}+\epsilon\right)^{2}+n\left(\frac{1}{2 \sqrt{n}}+2 \epsilon\right)^{2}\leq \frac{1}{2}+\frac{1}{2\sqrt{n}}+\frac{1}{2n}\leq \frac{7}{8}.
\] 

Then $\ds\frac{u}{(t-m)^{2}+u^{2}} \geq \frac{8}{7} u$.

\medskip

Now we are ready to outline the modified algorithm. We proceed as follows:

\begin{enumerate}
	\item Check if we have $k\geq n/2$ roots in a circle of small radius $\ds\epsilon$. If that is the case, find the center $c$ and radius $r$.
	
	\item If $k>n/2$, then from Proposition \ref{prop:3.1} we must have $u\leq 2\epsilon\sqrt{n}$ small and from Theorem \ref{thm:bound_more} we have 
	\[
	|t-c|\leq (2n+3)\epsilon.
	\] 	
	\item If $k=n/2$, call $r_1:=r$, $c_1:=c$ and construct the center $c_2$ and radius $r_2$ of the other $n/2$ roots.
	\begin{enumerate}
		\item if $r_2\leq \sqrt{\epsilon}$ and $|c_1-c_2|\leq r_2+r_1+\sqrt{\epsilon}+\epsilon$, then $|c_1-c_2|\leq 4\sqrt{\epsilon}$ and all the roots are on a circle centered at $c_1$ of radius $5\sqrt{\epsilon}$. Then from Theorem \ref{thm:bound_more}, we have:
		\[
		|t-c_1|\leq 5(2n+3)\sqrt{\epsilon}.
		\]
Note that we also have $u\leq 10\sqrt{n\epsilon}$ from Proposition \ref{prop:3.1}.
		
		\item if $r_2\leq \sqrt{\epsilon}$ and $|c_1-c_2|\geq r_2+r_1+\sqrt{\epsilon}+\epsilon$, we compute $\ds\frac{r_1}{r_2}.$
		
		\begin{itemize}
			\item If $\ds\frac{r_1}{r_2}\leq\frac{10n\epsilon}{3}$, we have from Theorem \ref{thm:u_large} that .
		
		\[
		|t-c_1|\leq \frac{1}{2\sqrt{n}}+\epsilon
		\]
		
Note that also $\ds u\leq \frac{20n\epsilon}{7}$ from Proposition \ref{prop:r2small_u}.
		
		\item  If $\ds\frac{r_1}{r_2}>\frac{10n\epsilon}{3}$, then $u$ is large enough to proceed with the usual algorithm, as Theorem \ref{thm:usmall_quot} implies $u>\epsilon$. 
		\end{itemize}
		
		\item if $r_2>\sqrt{\epsilon}$ and $|c_1-c_2|\geq 2r_2$, check if $\ds\frac{r_1}{r_2}<\frac{10n\epsilon}{3}$. If that is the case, from Proposition \ref{prop:r2notsmall_u} we have $u<8\sqrt{n}\epsilon$ and from Proposition \ref{prop:4.11} we have
		\[
		d_1\leq \frac{1}{2\sqrt{n}}+\epsilon.
		\]

		\item if $r_2>\sqrt{\epsilon}$ and $|c_1-c_2|\leq 2r_2$:
		
		\begin{enumerate}
		
		\item Check if we have more than $n/2$ roots on a circle of radius $2\epsilon$, with some center $c$. If that is the case, we have from Proposition \ref{prop:3.1} that $u\leq 4\sqrt{n}\epsilon$ and from Theorem \ref{thm:bound_more}: 
		\[
		|t-c|\leq 2(2n+3)\epsilon.
		\]
	If this is not that case:
		\item  if $c_1, c_2\not\in\RR$ and $c_1\neq \overline{c_2}$, or $c_1= \overline{c_2}$ and $r_1\neq r_2$, then we have $u>\epsilon$ from Proposition \ref{prop:4.1} and we can proceed with the usual algorithm.
		
		\item if $c_1=\overline{c_2}$ and $r_1=r_2$, then $2\text{Im}(c_1)=|c_1-c_2|\leq 2r_2\leq 2\epsilon$, thus all roots are on a circle of radius $2\epsilon$ centered at $\text{Re}(c_1)$, same as Case (i) above.

			\item  if $c_1, c_2\in \RR$ and $r_1r_2<\frac{3}{64n^2}$, then from Proposition \ref{prop:4.13} we have:
			
	\[
	|t-c_1|\leq \frac{1}{2\sqrt{n}}+\epsilon.
	\]
		
		\item If $c_1, c_2\in \RR$ and $\ds r_1r_2>\frac{3}{64n^2}$, we need to check with the usual algorithm if $u$ is close to $0$ or not. If $u$ is close to $0$, then we can use Theorem \ref{thm:u_small} to get:
		\[
		|t-c_1|\leq \frac{1}{2\sqrt{n}}+\epsilon.
		\]
		
	\end{enumerate}		
\end{enumerate}

\item In all cases when we got $u$ small, we can use Lemma \ref{lem:u_increase}
and take 
\[
m=\lfloor \text{Re}(c)+1/2\rfloor
\]
in the algorithm of \cite{SC}. We get a new binary form 
\[
F'(X, Z):=F(Z, -X-mZ)
\] 
with covariant point $z(F')=(t', u')=\left( \frac{-(t-m)}{(t-m)^{2}+u^{2}}, \frac{u}{(t-m)^{2}+u^{2}}\right)$ and $u'\geq 8u/7$.
	
\end{enumerate}

 As after every execution the value of $u$ increases and we have the bound $r\geq u/\sqrt{n}$, the value of $r$ also increases and we eventually exit the loop. In this case $u$ is also large enough to use the usual algorithm.

\end{document}